\documentclass[11pt]{amsproc}
\usepackage{amsmath}
\usepackage{amssymb}
\usepackage{amsfonts}

\theoremstyle{plain}

\newtheorem{definition}{Definition}

\newtheorem{lemma}{Lemma}

\newtheorem{proposition}{Proposition}

\newtheorem{theorem}{Theorem}
\numberwithin{equation}{section}
\input{tcilatex}

\begin{document}
\title[Nodal solutions to quasilinear elliptic equations...]{Nodal solutions
to quasilinear elliptic equations on compact Riemannian manifolds}
\author{Mohammed Benalili}
\address{Dept. of Mathematics BP119, Faculty of Sciences, University
Abou-Bekr Belka\"{\i}d Tlemcen Algeria}
\email{m\_benalili@yahoo.fr}

\begin{abstract}
We show the existence of nodal solutions to perturbed quasilinear elliptic
equations with critical Sobolev exponent on compact Riemannian manifolds. A
nonexistence result is also given.
\end{abstract}

\maketitle

\section{Introduction}

In this paper we investigate nodal solutions to quasilinear elliptic
equations involving terms with critical growth on compact manifolds. Nodal
solutions to scalar curvature type equation has been the subject of
investigation by various authors. Among them, we cite D. Holcman\cite{8}, A.
Jourdain \cite{9}, Z. Djadli and A. Jourdain\cite{5}. This work is an
extension to a previous one by Z. Djadli and A. Jourdain\cite{5} where the
authors studied the case of the Laplacian. We use variational methods based
on the Mountain Pass Theorem as done in H. Brezis and L. Nirenberg \cite{3}
and some ideas due to H. Hebey regarding isometry concentration. We approach
the problem via subcritical exponents, an idea originated by Yamabe. A non
existence result of nodal solution based on a Pohozahev type identity is
also given. Let $(M,g)$ be a smooth compact Riemannian manifold $n\geq 3$,
with or without boundary $\partial M$ and $p\in (1,n)$. We use the notations
of \cite{5}, let 
\begin{equation*}
W^{1,p}(M)=\left\{ 
\begin{array}{c}
H_{1}^{p}(M)\text{ \ if \ \ }\partial M=\phi \\ 
\overset{o}{H_{1}^{p}}(M)\text{ if \ }\partial M\neq \phi%
\end{array}%
\right.
\end{equation*}%
where $H_{1}^{p}(M)$ is the completion of $C^{\infty }(M)$ with respect to
the norm 
\begin{equation*}
\left\Vert u\right\Vert _{1,p}=\left\Vert \nabla u\right\Vert
_{p}+\left\Vert u\right\Vert _{p}
\end{equation*}%
and $\overset{o}{H_{1}^{p}}(M)$ is the completion of $C_{o}^{\infty }(M)$
with respect to the same norm. Let $G$ be a subgroup of the isometry group
of $(M,g)$ denoted $Isom(M)$. We assume that $G$ is compact. We also
consider $\tau $ an involutive isometry of $(M,g)$ that is an element of $%
Isom(M)$ such that $\tau o\tau =id_{M}$. For $x$ a point of $M$, we denote
by $O_{G}(x)$ the orbit of $x$ under the action of $G$. We say that $G$ and $%
\tau $ commute weakly if for every $x\in M$, $\tau (O_{G}(x))=O_{G}(\tau
(x)) $. We also say that the fixed points of $\tau $ splits $M$ into two
domains $\Omega _{1}$ and $\Omega _{2}$ stable under the action of $G$ if

(i) $\ M=\Omega _{1}\cup \Omega _{2}\cup \digamma $, with $\Omega _{1}\cap
\Omega _{2}=\phi $ and $mes(\digamma )=0.$

(ii) $\tau (\Omega _{1})=\Omega _{2}$ , and $\forall \sigma \in G$, $\forall 
$ $i=1,2$ \ $\sigma (\Omega _{i})=\Omega _{i}$ \ 

where $\digamma $ denotes the set of the fixed points of $\tau $, that is $%
\digamma =\left\{ x\in M:\tau (x)=x\right\} .$ We say that a function $u\in
W^{1,p}(M)$ is $\tau $- antisymmetric if $uo\tau =-u$ \ a.e and $G$%
-invariant if for all $\ \sigma \in G$, $uo\sigma =u$ a.e. In what follows,
we denote by $Card$ the cardinality of a set. We say that an operator $L$
defined on $W^{1,p}(M)$ is coercive on a subspace $X$ of $W^{1,p}(M)$ if
there exists a positive real $\Lambda $ such that for all $u\in X$, 
\begin{equation*}
\int_{M}L(u)udv_{g}\geq \Lambda \left\Vert u\right\Vert _{1,p}^{p}\text{.}
\end{equation*}%
Let $a$, $f$, $h$ be smooth functions on $M$, and $p^{\ast }=\frac{np}{n-p}$%
, $q\in (p-1,p^{\ast }-1)$, we consider the following equation%
\begin{equation}
\Delta _{p}u+a\left\vert u\right\vert ^{p-2}u=f\left\vert u\right\vert
^{p^{\ast }-2}u+h\left\vert u\right\vert ^{q-1}u  \tag{1}  \label{1}
\end{equation}%
with in case $M$ \ has a boundary $u=0$ on $\partial M$ ; where $\Delta
_{p}u=-div_{g}(\left\vert \nabla u\right\vert ^{p-2}\nabla u)$. Under
assumptions which will be precise later, we investigate nodal solutions of
equation(\ref{1}). By definition, a function $u\in W^{1,p}(M)$ is said to be
a weak solution of the equation(\ref{1}) if $u$ satisfies (\ref{1}) in the
distribution sense.

We say that the functions $a$, $f$, $h$ satisfy the conditions $(C)$ at an
interior point $x_{o}$ of $M$ if

\begin{equation*}
\left\{ 
\begin{array}{c}
(i)\ \ \ \ \ \ \ \ \ \ \ \ \ \ \ \ \ \ \ \ \ \ \ \ \ \ \ \ \ \ \ \ \ \ \
1<p<2\ \ \ \ \ \ \ \ a(x_{o})<0 \\ 
(ii)\ \ \ p=2\ \ \ \ \ \ \frac{4(n-1)}{n-2}a(x_{o})-Scal(x_{o})+(n-4)\frac{%
\Delta f(x_{o})}{f(x_{o})}\ \ <0 \\ 
(iii)\ \ \ 2\ \ <p<\frac{n}{2}\ \ \ \ \ \ \ \ \ \ \ \ \ \ \ \ \frac{\Delta
f(x_{o})}{f(x_{o})}<\frac{p}{n-3p+2}scal(x_{o}). \\ 
(iv)\ \ \text{For all }1<p<n\text{, \ }h(x_{o})=0\text{ \ and }\Delta
h(x_{o})\leq 0\text{.}%
\end{array}%
\right.
\end{equation*}

We set $N=p^{\ast }-1$ and let $q\in (p-1,N).$Our main result in this paper
reads as

\begin{theorem}
\label{th3} Let $G$ be a compact subgroup of the isometry group of $(M,g)$, $%
n\geq 3$, $\tau $ an involutive isometry of $(M,g)$ such that $G$ and $\tau $
commute weakly. Let $a$, $f$ and $h$ be three smooth $G$-invariant and $\tau 
$-invariant functions.

We assume that:

(1) The operator $\varphi \rightarrow $ $\Delta _{p}\varphi +a\left\vert
\varphi \right\vert ^{p-2}\varphi $ is coercive on the space $H=\left\{ u\in
W^{1,p}(M):\ u\text{ is }G\text{-invariant and }\tau \text{-antisymetric}%
\right\} $

(2) $f$ is positive on $M$ and attains its maximum at an interior point $%
x_{o}$ such that $\tau (O_{G}(x_{o}))\cap O_{G}(x_{o})=\phi $

(3) The functions $a$, $f$, $h$ satisfy the condition(C) at $x_{o}$.

Then equation(\ref{1}) possesses a nodal solution $u\in C^{1,\alpha }(M)$
which is $G$-invariant and $\tau $-antisymmetric. Moreover, if we assume
that the set $\digamma $ of fixed points of $\tau $ splits $M$ into two
domains $\Omega _{1}$and $\Omega _{2}$ stable under the action of$\ G$, we
can choose $u$ such that the zero set of $u$ is exactly $\digamma \cup
\partial M$.
\end{theorem}

\section{A generic theorem of existence}

First we give a regularity and a strong maximum results adapted to the
context of manifolds from those of Tolksdorf \cite{14}, Guedda-Veron \cite{7}
and Vasquez \cite{15} when dealing with Euclidian context. These results are
also given by Druet \cite{6} in the context of compact manifolds without
boundary. The proofs are similar and based on Moser's iteration scheme.

\begin{theorem}
\label{th1}($C^{1,\alpha }-$regularity) \ Let $\left( M,g\right) $ be a
compact Riemannian $n$- manifold with or without boundary, $n\geq 2$, $p\in
(1,n)$.

If $\ u\in W^{1,p}\left( M\right) $ is a solution of equation (\ref{1}) then 
$u\in C^{1,\alpha }\left( M\right) $.
\end{theorem}

\begin{proof}
Put 
\begin{equation}
g(x,u)=-a(x)\left\vert u\right\vert ^{p-2}u+f(x)\left\vert u\right\vert
^{p^{\ast }-2}u+h(x)\left\vert u\right\vert ^{q-1}u  \tag{2}  \label{2}
\end{equation}%
and 
\begin{equation*}
\widetilde{h}(x)=\frac{g(x,u(x))}{1+\left\vert u(x)\right\vert ^{p-1}}\text{.%
}
\end{equation*}%
Then%
\begin{equation*}
\left\vert \widetilde{h}(x)\right\vert \leq \left\Vert a\right\Vert _{\infty
}+\left\Vert f\right\Vert _{\infty }\left\vert u\right\vert ^{p^{\ast
}-p}+\left\Vert h\right\Vert _{\infty }\left\vert u\right\vert ^{q+1-p}
\end{equation*}%
where $\left\Vert .\right\Vert _{\infty }$ denotes the supremum norm. Since $%
u\in W^{1,p}(M)$, we have $\widetilde{h}\in L^{\frac{n}{p}}(M)$. The
equation(\ref{1}) reads as follows%
\begin{equation}
\Delta _{p}u=\left( 1+\left\vert u(x)\right\vert ^{p-1}\right) \widetilde{h}%
\text{.}  \tag{3}  \label{3}
\end{equation}%
Following arguments as in \ Guedda-Veron \cite{7} and Vasquez \cite{15} when
dealing with Euclidian context we first show that any solution $u\in
W^{1,p}(M)$ belongs to \ $L^{q}(M)$ for every $q\in \left[ 1,\infty \right[ $
. Let $k\geq 0$ and $v=\inf (\left\vert u\right\vert ,C)$ where $C$ is some
positive constant.

Multiplying equation(\ref{3}) by $v^{kp+1}$ and integrating over $M$, we get%
\begin{equation}
(kp+1)\int_{M}\left\vert u\right\vert ^{kp}\left\vert \nabla u\right\vert
^{p}dv_{g}=\int_{M}sgn(u)\widetilde{h}\left( 1+\left\vert u(x)\right\vert
^{p-1}\right) v^{kp+1}dv_{g}\text{.}  \tag{4}  \label{4}
\end{equation}%
On other hand, we have%
\begin{equation*}
\left\vert \nabla \left\vert u\right\vert ^{k+1}\right\vert
^{p}=(k+1)^{p}\left\vert u\right\vert ^{kp}\left\vert \nabla u\right\vert
^{p}
\end{equation*}%
so the equality (\ref{4}) writes%
\begin{equation}
\frac{kp+1}{(k+1)^{p}}\int_{M}\left\vert \nabla \left\vert u\right\vert
^{k+1}\right\vert ^{p}dv_{g}=\int_{M}sgn(u)\widetilde{h}\left( 1+\left\vert
u(x)\right\vert ^{p-1}\right) v^{kp+1}dv_{g}\text{.}  \tag{5}  \label{5}
\end{equation}%
Using Sobolev's inequality, we obtain for any fixed $\epsilon >0$%
\begin{equation*}
\left\Vert \left\vert u\right\vert ^{(k+1)}\right\Vert _{p^{\ast
}}^{p}=\left\Vert u\right\Vert _{(k+1)p^{\ast }}^{(k+1)p}
\end{equation*}%
\begin{equation*}
\leq \left( K(n,p)^{p}+\epsilon \right) \left\Vert \nabla \left\vert
u\right\vert ^{k+1}\right\Vert _{p}^{p}+B\left\Vert u\right\Vert
_{(k+1)p}^{(k+1)p}\text{.}
\end{equation*}%
Taking into account the relation(\ref{5}) and the H\"{o}lder's inequality we
get%
\begin{equation*}
\left\Vert u\right\Vert _{(k+1)p^{\ast }}^{(k+1)p}\leq \left(
K(n,p)^{p}+\epsilon \right) \frac{(k+1)^{p}}{kp+1}C^{kp+1}\left( \left\Vert
u\right\Vert _{(p-1)\frac{p^{\ast }}{p}}^{p-1}+Vol(M)^{1-\frac{p}{n}}\right)
\left\Vert \widetilde{h}\right\Vert _{\frac{n}{p}}
\end{equation*}%
\begin{equation*}
+B\left\Vert u\right\Vert _{(k+1)p}^{(k+1)p}
\end{equation*}%
where $K(n,p)$ is the best constant in the Sobolev's embedding $%
H_{1}^{p}(R^{n})\subset L^{p^{\ast }}(M)$ (see T. Aubin \cite{1} ) and $B$
is a positive constant depending on $\epsilon $.

Now, taking $(k+1)p=p^{\ast }$ i.e. $k=$ $\frac{p}{n-p}$, we obtain by the H%
\"{o}lder's inequality that%
\begin{equation*}
\left\Vert u\right\Vert _{p^{\ast }(1+\frac{p}{n-p})}\leq
\end{equation*}%
\begin{equation*}
\left\{ \left( K(n,p)^{p}+\epsilon \right) \frac{(k+1)^{p}}{kp+1}%
C^{kp+1}\left( Vol(M)^{\frac{1}{p^{\ast }}}+Vol(M)^{1-\frac{p}{n}}\right)
\left\Vert \widetilde{h}\right\Vert _{\frac{n}{p}}+B\right\} ^{\frac{1}{%
p\ast }}
\end{equation*}%
\begin{equation*}
\times \max (1,\left\Vert u\right\Vert _{p^{\ast }})\text{.}
\end{equation*}%
Consequently by a bootstrap arguments we get 
\begin{equation*}
u\in \dbigcap\limits_{1\leq q<\infty }L^{q}(M)\text{.}
\end{equation*}%
Now using the Moser's iteration scheme we are going to show that $u\in
L^{\infty }(M)$.

With the function $g$ given as in (\ref{2}), equation(\ref{1}) reads 
\begin{equation}
\Delta _{p}u=g\text{.}  \tag{6}  \label{6}
\end{equation}%
For any $k>1$, letting $t=k+p-1$, we get $\ \ \ $%
\begin{equation*}
\left\Vert \left\vert u\right\vert ^{\frac{t}{p}-1}\nabla u\right\Vert
_{p}^{p}=\int_{M}\left\vert t\right\vert ^{t-p}\left\vert \nabla
u\right\vert ^{p}dv_{g}\text{.}
\end{equation*}%
and multiplying equation(\ref{6}) by $\left\vert u\right\vert ^{k}$and
integrating over $M$, we obtain

\begin{equation*}
\int_{M}\left\vert u\right\vert ^{k}\Delta _{p}udv_{g}=k\int_{M}\left\vert
\nabla u\right\vert ^{p}\left\vert u\right\vert ^{k-2}udv_{g}
\end{equation*}%
\begin{equation}
=\int_{M}g\left\vert u\right\vert ^{k}dv_{g}\text{.}  \tag{7}  \label{7}
\end{equation}%
Using Sobolev's inequality, we obtain for any fixed $\epsilon >0$%
\begin{equation*}
\left\Vert \left\vert u\right\vert ^{\frac{t}{p}}\right\Vert _{p^{\ast
}}^{p}=\left\Vert u\right\Vert _{t\frac{p^{\ast }}{p}}^{t}
\end{equation*}%
\begin{equation*}
=\left( K(n,p)^{p}+\epsilon \right) \left( \frac{t}{p}\right) \left\Vert
\left\vert u\right\vert ^{\frac{t}{p}-1}\nabla u\right\Vert
_{p}^{p}+B\left\Vert u\right\Vert _{t}^{t}\text{.}
\end{equation*}

and since 
\begin{equation*}
\left\Vert \left\vert u\right\vert ^{\frac{t}{p}-1}\nabla u\right\Vert
_{p}^{p}=\int_{M}\left\vert t\right\vert ^{t-p}\left\vert \nabla
u\right\vert ^{p}dv_{g}
\end{equation*}%
and taking account of (\ref{7}) we obtain%
\begin{equation*}
\left\vert \int_{M}\left\vert u\right\vert ^{k}\Delta _{p}udv_{g}\right\vert
=k\int_{M}\left\vert u\right\vert ^{k-1}\left\vert \nabla u\right\vert
^{p}dv_{g}=k\int_{M}\left\vert u\right\vert ^{t-p}\left\vert \nabla
u\right\vert ^{p}dv_{g}
\end{equation*}%
\begin{equation*}
=k\left\Vert \left\vert u\right\vert ^{\frac{t}{p}-1}\nabla u\right\Vert
_{p}^{p}\leq \left\Vert g\right\Vert _{s}\left\Vert u\right\Vert _{kr}^{k}
\end{equation*}%
where $r,s>1$ are conjugate numbers.
\end{proof}

Consequently%
\begin{equation*}
\left\Vert u\right\Vert _{t\frac{p^{\ast }}{p}}^{t}\leq \left(
K(n,p)^{p}+\epsilon \right) \left( \frac{t}{p}\right) ^{p}\left\Vert
g\right\Vert _{s}\left\Vert u\right\Vert _{kr}^{k}+B\left\Vert u\right\Vert
_{t}^{t}
\end{equation*}%
and by H\"{o}lder's inequality we get%
\begin{equation*}
\left\Vert u\right\Vert _{t\frac{p^{\ast }}{p}}^{t}\leq \left(
K(n,p)^{p}+\epsilon \right) \left( \frac{t}{p}\right) ^{p}\left\Vert
g\right\Vert _{s}\left\Vert u\right\Vert _{rt}^{k}Vol(M)^{\frac{p-1}{rt}}
\end{equation*}%
\begin{equation*}
+B\left\Vert u\right\Vert _{rt}^{t}Vol(M)^{1-\frac{1}{r}}\text{.}
\end{equation*}%
Then%
\begin{equation*}
\left\Vert u\right\Vert _{t\frac{p^{\ast }}{p}}\leq \left( \frac{t}{p}%
\right) ^{\frac{^{p}}{t}}Vol(M)^{\frac{p-1}{rt}}\max \left\{ \left(
K(n,p)^{p}+\epsilon \right) ,B\right\}
\end{equation*}%
\begin{equation*}
\times \left( \left\Vert g\right\Vert _{s}+Vol(M)^{^{\frac{(r-1)t-p+1}{rt}%
}}\right) ^{\frac{1}{t}}\max (1,\left\Vert u\right\Vert _{rt}^{t})
\end{equation*}%
or%
\begin{equation}
\left\Vert u\right\Vert _{t\frac{p^{\ast }}{p}}\leq \left( \frac{t}{p}%
\right) ^{\frac{^{p}}{t}}A^{\frac{1}{t}}\text{ }\max \left( 1,\left\Vert
u\right\Vert _{rt}\right)  \tag{8}  \label{8}
\end{equation}%
where $A$ \ is a constant independent of $t$. Now we choose $r<\frac{p^{\ast
}}{p}=\frac{n}{n-p}$ i.e. $s>\frac{n}{p}$ which is possible by the first
part of the proof.

\begin{proof}
Let $\alpha >0$ such that $r(1+\alpha )=\frac{P^{\ast }}{p}$and $\beta
=1+\alpha $. Let also $t=\beta ^{i}$ where $i$ is a positive integer; the
recurrent \ relation(\ref{8}) writes as%
\begin{equation*}
\left\Vert u\right\Vert _{r\beta ^{i+1}}\leq \left( \frac{\beta ^{i}}{p}%
\right) ^{\frac{^{p}}{\beta ^{i}}}A^{\frac{1}{\beta ^{i}}}\max \left(
1,\left\Vert u\right\Vert _{r\beta ^{i}}\right)
\end{equation*}%
and recurrently we get%
\begin{equation*}
\left\Vert u\right\Vert _{r\beta ^{i+1}}\leq A^{\sum_{j=1}^{i}\frac{1}{\beta
^{j}}}\frac{\beta ^{p(\sum_{j=1}^{i}\frac{j}{\beta ^{j}})}}{p(\sum_{j=1}^{i}%
\frac{1}{\beta ^{j}})}\max \left( 1,\left\Vert u\right\Vert _{r}\right) 
\text{.}
\end{equation*}%
Now since the series $\sum_{i=1}^{\infty }\frac{1}{\beta ^{i}}=\frac{1}{%
\beta -1}$ and \ $\sum_{j=1}^{\infty }\frac{j}{\beta ^{j}}$ are convergent,
we get \ that $u\in L^{\infty }(M)$. At this stage the conclusion follows
from theorem of P. Tolksdorf\cite{14}
\end{proof}

\begin{theorem}
\label{th2}(Strong maximum principle) Let $\left( M,g\right) $ be a compact
Riemannian $n$-manifold with or without boundary, $p\in (1,n),$ and let $%
u\in C_{o}^{1}\left( M\right) $ be such that 
\begin{equation*}
\Delta _{p}u+f\left( .,u\right) \geq 0\text{ \ \ on }M\text{,}
\end{equation*}%
$f$ such that 
\begin{equation*}
\left\{ 
\begin{array}{c}
f\left( x,r\right) <f\left( x,s\right) ,\text{ \ \ \ }\forall x\in M\text{\
\ \ }\forall 0\leq r<s \\ 
\left\vert f\left( x,r\right) \right\vert \leq C\left( K+\left\vert
r\right\vert ^{p-2}\right) \left\vert r\right\vert ,\text{ \ \ \ }\forall
\left( x,r\right) \in M\times R%
\end{array}%
\right.
\end{equation*}%
where $C$ and $K$ are positive constants.

If $u\geq 0$ on $M$ and $u$ does not vanish identically, then $u>0$ on $\
int(M)=M-\partial M$.$.$
\end{theorem}

Let $G$ be a compact subgroup of the isometry group of $(M,\,g)$ and $\tau $
be an involutive isometry of $(M,g)$. We assume that $G$ and $\tau $ commute
weakly for some $x_{1}\in M$, $\tau (O_{G}(x_{1}))\cap O_{G}(x_{1})=\phi $.
Then 
\begin{equation*}
H=\left\{ u\in W^{1,p}(M)\text{, }u\text{ is }G\text{-invariant and }\tau -%
\text{antisymmetric }\right\}
\end{equation*}%
is not trivial. Indeed $H$ contains the test function given in section3.

Denote by $\left\langle G,\tau \right\rangle $ the subgroup of the isometry
group $Isom(M,g)$ generated by $G$ and $\tau $ and by $K(n,p)$ the best
constant in the Sobolev's embedding of $W^{1,p}(R^{n})$ in $L^{\frac{pn}{n-p}%
}(R^{n})$.

We consider the following functional $J$ defined on \ $H$ \ \ \ by 
\begin{equation*}
J(\varphi )=\int_{M}\left\{ \frac{1}{p}\left\vert \nabla \varphi \right\vert
^{p}+a\frac{1}{p}\left\vert \varphi \right\vert ^{p}-\frac{1}{p^{\ast }}%
f\left\vert \varphi \right\vert ^{p^{\ast }}-\frac{1}{q+1}h\left\vert
\varphi \right\vert ^{q+1}\right\} dv_{g}\text{.}
\end{equation*}%
In this section we establish the following generic theorem.

\begin{theorem}
\label{th4} Let $G$ be a compact subgroup of the isometry group of $(M,g)$, $%
n\geq 3$, $\tau $ an involutive isometry of $(M,g)$ such that $G$ and $\tau $
commute weakly and $\tau (O_{G}(x_{1}))\cap O_{G}(x_{1})=\phi $ for some $%
x_{1}\in M$. Let also $a$, $f$ and $h$ be three smooth $G$-invariant and $%
\tau $-invariant \ functions. We assume that $f$ is positive on $M$ and the
operator $\varphi \rightarrow $ $\Delta _{p}\varphi +a\left\vert \varphi
\right\vert ^{p-2}\varphi $ is coercive on $H$. We set $N=p^{\ast }-1$ and
let $q\in (p-1,N).$We assume that for all $x$ in $M$ there exists $v\in H$, $%
v\neq 0$ such that 
\begin{equation}
\sup_{t\geq 0}\left\{ J(tv)\right\} <\frac{CardO_{\left\langle G,\tau
\right\rangle }(x)}{K(n,p)^{n}f(x)^{\frac{n-p}{p}}}\text{.}  \tag{8'}
\label{8'}
\end{equation}%
Then equation(\ref{1}) possesses a nodal solution $u\in C^{1,\alpha }(M)$
which is $G$-invariant and $\tau $-antisymmetric. Moreover, if we assume
that the set $\digamma $ of fixed points of $\tau $ splits $M$ into two
domains $\Omega _{1}$and $\Omega _{2}$ stable under the action of $\ G$, we
can choose $u$ such that the zero set of $u$ is exactly $\digamma \cup
\partial M$.
\end{theorem}

\subsection{\ The subcritical case\ }

Now, following the strategy originated by Yamabe, we prove the existence of
a nodal solution to the equation(\ref{1}) for the subcritical exponent.

\begin{proposition}
\label{p1} Let \ $G$ be a compact subgroup of the isometry group of ($M,g)$, 
$n\geq 3$, let $\tau $ be an involutive isometry of ($M,g$) such that $G$
and $\tau $ commute weakly and such that for some $x_{1}\in M$ \ $\tau
(O_{G}(x_{1}))\cap O_{G}(x_{1})=\phi $. Let also $a$, $f$ and $h$ be three
smooth $G$-invariant and $\tau $-invariant \ functions. We assume that $f$
is positive on $M$ and that the operator $\varphi \rightarrow \Delta
_{p}\varphi +a$ $\left\vert \varphi \right\vert ^{p-2}\varphi $ is coercive
on $H$. We set $N=p^{\ast }-1$ , $q\in (p-1,N)$ and let $\epsilon _{o}$ be
such that $0<\epsilon _{o}\leq N-q$. Then for all $\epsilon $ such that $%
0<\epsilon \leq \epsilon _{o}$ there exists $\varphi _{\epsilon }\in
C^{1,\alpha }(M)$, $G$-invariant and $\tau $-antisymmetric $\varphi
_{\epsilon }\neq 0$ in $M$ and $\varphi _{\epsilon }=0$ on $\partial M$
which is a nodal weak solution of the equation 
\begin{equation}
\Delta _{p}\varphi _{\epsilon }+a\left\vert \varphi _{\epsilon }\right\vert
^{p-2}u=f\left\vert \varphi _{\epsilon }\right\vert ^{p^{\ast }-2-\epsilon
}\varphi _{\epsilon }+h\left\vert \varphi _{\epsilon }\right\vert
^{q-2}\varphi _{\epsilon }\text{.}  \tag{9}  \label{9}
\end{equation}%
Moreover, if we assume that the set$\ \digamma $ of fixed points of $\tau $
splits $M$ into two domains $\Omega _{1}$ and $\Omega _{2}$ stable under the
action of $G$, we can choose $\varphi _{\epsilon }$ such that its zero set
is exactly $\digamma \cup \partial M$.
\end{proposition}

The proof of the Proposition(\ref{p1}) relies on the following Mountain-Pass
Lemma of Ambrosetti and Rabinowitz(\cite{2})

\begin{lemma}
Let $\phi $ be a $C^{1}$ function on a Banach space $X$. Suppose that there
exists a neighborhood $U$ of $0$ in $X$ , a $v\in X\backslash U$ and a
constant $\rho $ such that $\phi (0)<\rho $, $\phi (v)<\rho $ and $\phi
(u)\geq \rho $ for all $u\in \partial U$. Let $\Gamma $ denote the class of
continuous paths joining $0$ to $v$ and $c=\inf_{\gamma \in \Gamma
}\max_{w\in \gamma }\phi (w)$.

Then there is a sequence $\left( u_{j}\right) _{j}$ in $X$ such that $\phi
(u_{j})\rightarrow c$ and $\phi ^{\prime }(u_{j})\rightarrow 0$ in $X^{\ast
} $ ( dual space of $X$).
\end{lemma}

We recall the following concepts

\begin{definition}
Let $X$ be a Banach and $\phi $ a function of class $C^{1}$ on $X$. We say
that $u_{n}\in X$ is a Palais-Smale sequence at level $c$ ( Shortly a (PS)$%
_{c}$ sequence ), if

(i) $\ \ \phi (u_{n})\rightarrow c$

(ii) $\phi ^{\prime }(u_{n})\rightarrow 0$.
\end{definition}

\begin{definition}
We say that $\phi $ satisfies the (PS)$_{c}$ condition if every (PS)$_{c}$
sequence has a converging subsequence . We say that $\phi $ satisfies the
(PS) condition if it satisfies the (PS)$_{c}$ condition for all $c\in R$.
\end{definition}

For any sufficiently small $\epsilon $ such that \ $0<\epsilon \leq \epsilon
_{o}$, we consider the $C^{1}$-functional $J_{\epsilon }$ defined on the
space $H$ by 
\begin{equation*}
J_{\epsilon }(\varphi )=\int_{M}\left\{ \frac{1}{p}\left\vert \nabla \varphi
\right\vert ^{p}+\frac{1}{p}a\left\vert \varphi \right\vert ^{p}-\frac{1}{%
p^{\ast }-\epsilon }f\left\vert \varphi \right\vert ^{p^{\ast }-\epsilon }-%
\frac{1}{q+1}h\left\vert \varphi \right\vert ^{q+1}\right\} dv_{g}\text{.}
\end{equation*}%
Following Brezis and Nirenberg (\cite{3}), we show that the functional $%
J_{\epsilon }$ satisfies the assumptions of the Mountain-Pass Lemma, for
every $\epsilon $ such that \ $0<\epsilon \leq \epsilon _{o}$.

\begin{lemma}
\label{lem1} For every $\epsilon $ such that \ $0<\epsilon \leq \epsilon
_{o} $, there exists a ball $U$ of radius independent of $\epsilon $ around $%
0$ in $H$ included in the unit ball, and a positive real number $\rho $
independent of $\ \epsilon $ such that

(i) $\ \ \ \ \ \ \ \ \forall $ $\varphi \in U$ , $J_{\epsilon }(\varphi
)\geq \rho >0$

(ii) \ \ \ \ \ \ \ \ $\exists \psi \notin U$ such that \ $J_{\epsilon }(\psi
)<\rho .$
\end{lemma}

\begin{proof}
By the coercivity of the operator $\varphi \rightarrow \Delta _{p}\varphi
+a\left\vert \varphi \right\vert ^{p-2}\varphi $, there exists a positive
real number $\Lambda $ such that 
\begin{equation*}
J_{\epsilon }(\varphi )\geq \frac{\Lambda }{p}\left\Vert \varphi \right\Vert
_{1,p}^{p}-\frac{1}{q+1}\left\Vert h\right\Vert _{\infty }\left\Vert \varphi
\right\Vert _{q+1}^{q+1}-\frac{1}{p^{\ast }-\epsilon }\left\Vert
f\right\Vert _{\infty }\left\Vert \varphi \right\Vert _{p^{\ast }-\epsilon
}^{p^{\ast }-\epsilon }
\end{equation*}%
and by the Sobolev's inequality, that is for every $\eta >0$ 
\begin{equation*}
\left\Vert \varphi \right\Vert _{q+1}^{q+1}\leq \left[ \left(
K(n,p)^{p}+\eta \right) \left\Vert \nabla \varphi \right\Vert
_{p}^{p}+B\left\Vert \varphi \right\Vert _{p}^{p}\right] ^{\frac{q+1}{p}}%
\text{,}
\end{equation*}%
one has%
\begin{equation*}
J_{\epsilon }(\varphi )\geq \left\Vert \varphi \right\Vert _{1,p}^{p}\left[ 
\frac{\Lambda }{p}-\frac{1}{q+1}\left\Vert h\right\Vert _{\infty }\max
(\left( K(n,p)^{p}+\eta \right) ,B)^{\frac{q+1}{p}}\left\Vert \varphi
\right\Vert _{1,p}^{q+1-p}\right.
\end{equation*}%
\begin{equation*}
\left. -\frac{1}{p^{\ast }-\epsilon }\left\Vert f\right\Vert _{\infty }\max
(\left( K(n,p)^{p}+\eta \right) ,B)^{\frac{p^{\ast }-\epsilon }{p}%
}\left\Vert \varphi \right\Vert _{1,p}^{p^{\ast }-p-\epsilon }\right]
\end{equation*}%
and since $q+1-p>0$ and \ $p^{\ast }-p-\epsilon >0$, \ there is a ball $U$
included in the unit ball and a positive number $\rho $ independent of $%
\epsilon $ with $0<\epsilon \leq \epsilon _{o}$ such that for every $u\in
\partial U$, \ $J_{\epsilon }(\varphi )\geq \rho .$

For $t>0$, 
\begin{equation*}
J_{\epsilon }(t\varphi )\leq t^{p}\left\{ \frac{1}{p}\left\Vert \nabla
\varphi \right\Vert _{p}^{p}+\frac{1}{p}\left\Vert a\right\Vert _{\infty
}\left\Vert \varphi \right\Vert _{p}^{p}-\frac{t^{q+1-p}}{q+1}\min_{x\in
M}h(x)\left\Vert \varphi \right\Vert _{q+1}^{q+1}\right.
\end{equation*}%
\begin{equation*}
\left. -\frac{t^{p^{\ast }-p-\epsilon }}{p^{\ast }-\epsilon }\min_{x\in
M}f(x)\left\Vert \varphi \right\Vert _{p^{\ast }-\epsilon }^{p^{\ast
}-\epsilon }\right\}
\end{equation*}%
so since $p^{\ast }-q-1>0$, there is a sufficiently large $t_{o}$ such that
\ if $\psi =t_{o}\varphi $ then $\psi \notin U$ and $J_{\epsilon }(\psi
)<\rho $ for $\varepsilon $ sufficiently small.
\end{proof}

Let $P$ the class of continuous paths joining $0$ to $\psi $and let $%
c_{\epsilon }=\inf_{\gamma \in P}\max_{w\in \gamma }J_{\epsilon }(w)$. Then
by Lemma\ref{lem1} there exists a (PS)$_{c_{\epsilon }\text{ }}$sequence in $%
H$.

Now we are going to show that each $(PS)$ sequence satisfies the
Palais-Smale condition.

\begin{lemma}
\label{lem2} Each Palais-Smale sequence for the functional $J_{\epsilon }$
is bounded.
\end{lemma}

\begin{proof}
We argue by contradiction. Suppose that there exists a sequence $\left\{
\varphi _{j}\right\} $ such that $J_{\epsilon }(\varphi _{j})$ tends to a
finite limit $c,$ $J_{\epsilon }^{^{\prime }}(\varphi _{j})$ goes to zero
and $\varphi _{j}$ to infinite in the $W^{1,p}(M)$-norm. More explicitly we
have for each $\psi \in W^{1,p}(M)$%
\begin{equation*}
\int_{M}\left\{ \frac{1}{p}\left\vert \nabla \varphi _{j}\right\vert ^{p}+%
\frac{1}{p}a\left\vert \varphi _{j}\right\vert ^{p}-\frac{1}{p^{\ast
}-\epsilon }f\left\vert \varphi _{j}\right\vert ^{p^{\ast }-\epsilon }-\frac{%
1}{q+1}h\left\vert \varphi _{j}\right\vert ^{q+1}\right\} dv_{g}\rightarrow c
\end{equation*}%
and 
\begin{equation*}
\int_{M}\left\vert \nabla \varphi _{j}\right\vert ^{p-2}\nabla _{i}\varphi
_{j}\nabla ^{i}\psi dv_{g}+\int_{M}a\left\vert \varphi _{j}\right\vert
^{p-2}\varphi _{j}\psi dv_{g}-\int_{M}f\left\vert \varphi _{j}\right\vert
^{p^{\ast }-2-\epsilon }\varphi _{j}\psi dv_{g}
\end{equation*}%
\begin{equation*}
-\int_{M}h\left\vert \varphi _{j}\right\vert ^{q-1}\varphi _{j}\psi
dv_{g}\rightarrow 0
\end{equation*}%
so for any $\eta >0$ there exists a positive integer $N$ such that for every 
$j\geq N$ one has 
\begin{equation*}
\left\vert \int_{M}\left\{ \frac{1}{p}\left\vert \nabla \varphi
_{j}\right\vert ^{p}+\frac{1}{p}a\left\vert \varphi _{j}\right\vert ^{p}-%
\frac{1}{p^{\ast }-\epsilon }f\left\vert \varphi _{j}\right\vert ^{p^{\ast
}-\epsilon }-\frac{1}{q+1}h\left\vert \varphi _{j}\right\vert ^{q+1}\right\}
dv_{g}-c\right\vert \leq \eta
\end{equation*}%
and 
\begin{equation*}
\left\vert \int_{M}\left\vert \nabla \varphi _{j}\right\vert ^{p-2}\nabla
_{i}\varphi _{j}\nabla ^{i}\psi dv_{g}+\int_{M}a\left\vert \varphi
_{j}\right\vert ^{p-2}\varphi _{j}\psi dv_{g}-\int_{M}f\left\vert \varphi
_{j}\right\vert ^{p^{\ast }-2-\epsilon }\varphi _{j}\psi dv_{g}\right.
\end{equation*}%
\begin{equation*}
\left. -\int_{M}h\left\vert \varphi _{j}\right\vert ^{q-1}\varphi _{j}\psi
dv_{g}\right\vert \leq \eta
\end{equation*}%
In the particular case where $\psi =\varphi _{j},$ we get 
\begin{equation*}
\left\vert \int_{M}\left\{ \frac{1}{p}\left\vert \nabla \varphi
_{j}\right\vert ^{p}+\frac{1}{p}a\left\vert \varphi _{j}\right\vert ^{p}-%
\frac{1}{p^{\ast }-\epsilon }f\left\vert \varphi _{j}\right\vert ^{p^{\ast
}-\epsilon }-\frac{1}{q+1}h\left\vert \varphi _{j}\right\vert ^{q+1}\right\}
dv_{g}-c\right\vert \leq \eta
\end{equation*}%
and 
\begin{equation*}
\left\vert \int_{M}\left\{ \left\vert \nabla \varphi _{j}\right\vert
^{p}dv_{g}+a\left\vert \varphi _{j}\right\vert ^{p}dv_{g}-f\left\vert
\varphi _{j}\right\vert ^{p^{\ast }-\epsilon }\right\} dv_{g}\right.
\end{equation*}%
\begin{equation}
\left. -\int_{M}h\left\vert \varphi _{j}\right\vert ^{q+1}dv_{g}\right\vert
\leq \eta \text{.}  \tag{10}  \label{10}
\end{equation}%
Then, we obtain 
\begin{equation}
\left\vert (1-\frac{p}{q+1})\int_{M}\left( \left\vert \nabla \varphi
_{j}\right\vert ^{p}+a\left\vert \varphi \right\vert _{j}^{p}\right)
dv_{g}+p(\frac{1}{q+1}-\frac{1}{p^{\ast }-\epsilon })\int_{M}f\left\vert
\varphi _{j}\right\vert ^{p^{\ast }-\epsilon }dv_{g}-pc\right\vert  \tag{11}
\label{11}
\end{equation}%
\begin{equation*}
\leq (1+\frac{1}{q+1})p\eta
\end{equation*}%
and 
\begin{equation}
\left\vert (1-\frac{p}{p^{\ast }-\epsilon })\int_{M}f\left\vert \varphi
_{j}\right\vert ^{p^{\ast }-\epsilon }dv_{g}+(1-\frac{p}{q+1}%
\int_{M}h\left\vert \varphi _{j}\right\vert ^{q+1}dv_{g}-pc\right\vert 
\tag{12}  \label{12}
\end{equation}%
\begin{equation*}
\leq (1+p)\eta \text{.}
\end{equation*}%
Now, since $f>0$, there is a constant $C>0$ such that 
\begin{equation*}
C(1-\frac{p}{p^{\ast }-\epsilon })\int_{M}\left\vert \varphi _{j}\right\vert
^{p^{\ast }-\epsilon }dv_{g}\leq (1-\frac{p}{q+1})\left\Vert h\right\Vert
_{\infty }\int_{M}\left\vert \varphi _{j}\right\vert
^{q+1}dv_{g}+pc+(1+p)\eta
\end{equation*}%
where $\left\Vert h\right\Vert _{\infty }=\max_{x\in M}\left\vert
h(x)\right\vert $.

On the other hand since $p^{\ast }-\epsilon >q+1$, for any $\nu >0$, there
exists a constant $C_{\nu }^{\prime }$ such that $t^{q+1}\leq \nu t^{p^{\ast
}-\epsilon }+C_{\nu }^{\prime }$ for any $t\geq 0$. So 
\begin{equation*}
C(1-\frac{p}{p^{\ast }-\epsilon })\int_{M}\left\vert \varphi _{j}\right\vert
^{p^{\ast }-\epsilon }dv_{g}\leq (1-\frac{p}{q+1})\left\Vert h\right\Vert
_{\infty }\left( \nu \int_{M}\left\vert \varphi _{j}\right\vert ^{p^{\ast
}-\epsilon }dv_{g}+C_{\nu }vol(M)\right)
\end{equation*}%
\begin{equation*}
+pc+(1+p)\eta
\end{equation*}%
and 
\begin{equation*}
\left[ C(1-\frac{p}{p^{\ast }-\epsilon })-(1-\frac{p}{q+1})\left\Vert
h\right\Vert _{\infty }\nu \right] \int_{M}\left\vert \varphi
_{j}\right\vert ^{p^{\ast }-\epsilon }dv_{g}\leq cste\text{.}
\end{equation*}%
Choosing $\nu >0$ small enough so that $C(1-\frac{p}{p^{\ast }-\epsilon }%
)-(1-\frac{p}{q+1})\left\Vert h\right\Vert _{\infty }\nu >0$ and get 
\begin{equation}
\int_{M}\left\vert \varphi _{j}\right\vert ^{p^{\ast }-\epsilon }dv_{g}\leq
cste\text{.}  \tag{13}  \label{13}
\end{equation}%
By Lemma \ref{lem1}, we can choose $\rho $ to be an $W^{1,p}(M)$- norm such
that 
\begin{equation*}
\inf_{\left\Vert \varphi \right\Vert _{1,p}=\rho }J_{\epsilon }(\varphi )>0%
\text{.}
\end{equation*}%
Letting $\psi _{j}=\rho \frac{\varphi _{j}}{\left\Vert \varphi
_{j}\right\Vert _{1,p}}$, we obtain from (\ref{13}) that

\begin{equation}
\int_{M}\left\vert \psi _{j}\right\vert ^{p^{\ast }-\epsilon }dv_{g}=O\left( 
\frac{\rho ^{p^{\ast }-\epsilon }}{\left\Vert \varphi _{j}\right\Vert
_{1,p}^{p^{\ast }-\epsilon }}\right)  \tag{14}  \label{14}
\end{equation}%
and by (\ref{11}), we get 
\begin{equation*}
\left\vert (1-\frac{p}{q+1})\int_{M}\left( \left\vert \nabla \psi
_{j}\right\vert ^{p}+a\left\vert \psi _{j}\right\vert ^{p}\right)
dv_{g}+\right.
\end{equation*}%
\begin{equation}
\left. p(\frac{1}{q+1}-\frac{1}{p^{\ast }-\epsilon })\frac{\left\Vert
\varphi _{j}\right\Vert _{1,p}^{p^{\ast }-\epsilon -p}}{\rho ^{p^{\ast
}-\epsilon }}\int_{M}f\left\vert \psi _{j}\right\vert ^{p^{\ast }-\epsilon
}dv_{g}-pc\right\vert  \tag{15}  \label{15}
\end{equation}%
\begin{equation*}
\leq (1+\frac{1}{q+1})p\eta
\end{equation*}%
Letting $j$ go to infinity, we obtain that $J_{\epsilon }(\psi _{j})$ tends
to zero. And since $\left\Vert \psi _{j}\right\Vert _{1,p}=\rho $, we have 
\begin{equation*}
\inf_{\left\Vert \varphi \right\Vert _{1,p}=\rho }J_{\epsilon }(\varphi
)\leq J_{\epsilon }(\psi _{j})
\end{equation*}%
so 
\begin{equation*}
\inf_{\left\Vert \varphi \right\Vert _{1,p}=\rho }J_{\epsilon }(\varphi
)\leq 0
\end{equation*}%
hence a contradiction. Then the sequence $\left\{ \varphi _{j}\right\} $ is
bounded in $W^{1,p}(M).$ Now since $q<p^{\ast }-1$, the Sobolev injections
are compact, so the Palais-Smale condition is satisfied. \ \ 
\end{proof}

Now we are in position to prove Proposition\ref{p1}.

\begin{proof}
( of Proposition\ref{p1} ) \ Let $C$ \ be the set of paths $\gamma $ joining 
$0$ and $\psi $ and $c_{\epsilon }=\inf_{\gamma \in C}\sup_{t\in \left[ 0,1%
\right] }\gamma (t).$ As a consequence of Lemma\ref{lem1}, there exists a
sequence $\left\{ \varphi _{j}\right\} \subset H$ such that $J_{\epsilon
}(\varphi _{j})\rightarrow c_{\epsilon }$ and $J_{\epsilon }^{\prime
}(\varphi _{j})\rightarrow 0$ strongly in $H^{\prime }.$ By Lemma\ref{lem2},
we can extract a subsequence still denoted $\left\{ \varphi _{j}\right\} $
such that

$\varphi _{j}\rightarrow \varphi _{\epsilon }$ weakly in $H$

$\varphi _{j}\rightarrow \varphi _{\epsilon }$ strongly in $L^{r\text{ \ \ \
\ }}$for $r$ $<p^{\ast }$

$\varphi _{j}\rightarrow \varphi _{\epsilon }$ a.e. in $M$.

Clearly $\varphi _{\epsilon }$ is a weak solution of the equation 
\begin{equation}
\Delta _{p}\varphi _{\epsilon }+a\left\vert \varphi _{\epsilon }\right\vert
^{p-2}\varphi _{\epsilon }-h\left\vert \varphi _{\epsilon }\right\vert
^{q-1}\varphi _{\epsilon }-f\left\vert \varphi _{\epsilon }\right\vert
^{p^{\ast }-2-\epsilon }\varphi _{\epsilon }=0\ \ \ \ \ \text{in }H\text{.} 
\tag{16}  \label{16}
\end{equation}

Now, to use the results of regularity, we must show that $\ \varphi
_{\epsilon }$ satisfies the equation(\ref{16}) weakly in $W^{1,p}(M)$. So we
consider $\psi \in W^{1,p}(M)$ and the Haar measure denoted by $d\sigma $ on
the isometric group $G$ . Set $\overset{\_}{\psi }(x)=\int_{G}\psi (\sigma
(x))d\sigma $ for all $x\in M$, then $\overset{\_}{\psi }$ is $G$- invariant
\ and it follows by multiplying the equation(\ref{16}) by $\overset{\_}{\psi 
}$ and integrating over $M$ that 
\begin{equation*}
\int_{M}\left\{ \left\vert \nabla \varphi _{\epsilon }\right\vert
^{p-2}\nabla _{i}\varphi _{\epsilon }\nabla ^{i}\overset{\_}{\psi }%
+a\left\vert \varphi _{\epsilon }\right\vert ^{p-2}\varphi _{\epsilon }%
\overset{\_}{\psi }-h\left\vert \varphi _{\epsilon }\right\vert
^{q-1}\varphi _{\epsilon }\overset{\_}{\psi }-f\left\vert \varphi _{\epsilon
}\right\vert ^{p^{\ast }-2-\epsilon }\varphi _{\epsilon }\overset{\_}{\psi }%
\right\} =0
\end{equation*}%
but 
\begin{equation*}
\nabla \int_{G}\psi (\sigma (x))d\sigma =\int_{G}\nabla \psi (\sigma
(x))d\sigma 
\end{equation*}%
then\ 
\begin{equation*}
\int_{M}\left\{ \left\vert \nabla \varphi _{\epsilon }\right\vert
^{p-2}\nabla _{i}\varphi _{\epsilon }\nabla ^{i}\left( \int_{G}\psi (\sigma
(x))d\sigma \right) +a\left\vert \varphi _{\epsilon }\right\vert
^{p-2}\varphi _{\epsilon }\left( \int_{G}\psi (\sigma (x))d\sigma \right)
\right. 
\end{equation*}%
\begin{equation*}
\left. -h\left\vert \varphi _{\epsilon }\right\vert ^{q-1}\varphi _{\epsilon
}\left( \int_{G}\psi (\sigma (x))d\sigma \right) -f\left\vert \varphi
_{\epsilon }\right\vert ^{p^{\ast }-2-\epsilon }\varphi _{\epsilon }\left(
\int_{G}\psi (\sigma (x))d\sigma \right) \right\} dv_{g}
\end{equation*}%
\begin{equation*}
=\int_{M}\int_{G}\left\{ \left\vert \nabla \varphi _{\epsilon }\right\vert
^{p-2}\nabla _{i}\varphi _{\epsilon }\nabla ^{i}\psi (\sigma
(x))+a\left\vert \varphi _{\epsilon }\right\vert ^{p-2}\varphi _{\epsilon
}\psi (\sigma (x))\right. 
\end{equation*}%
\begin{equation*}
\left. -h\left\vert \varphi _{\epsilon }\right\vert ^{q-1}\varphi _{\epsilon
}\psi (\sigma (x))-f\left\vert \varphi _{\epsilon }\right\vert ^{p^{\ast
}-2-\epsilon }\varphi _{\epsilon }\psi (\sigma (x))\right\} d\sigma dv_{g}=0%
\text{.}
\end{equation*}%
Now, by the Fubini's theorem, we get%
\begin{equation*}
\int_{G}\int_{M}\left\{ \left\vert \nabla \varphi _{\epsilon }\right\vert
^{p-2}\nabla _{i}\varphi _{\epsilon }\nabla ^{i}\psi (\sigma
(x))+a\left\vert \varphi _{\epsilon }\right\vert ^{p-2}\varphi _{\epsilon
}\psi (\sigma (x))\right. 
\end{equation*}%
\begin{equation*}
\left. -h\left\vert \varphi _{\epsilon }\right\vert ^{q-1}\varphi _{\epsilon
}\psi (\sigma (x))-f\left\vert \varphi _{\epsilon }\right\vert ^{p^{\ast
}-2-\epsilon }\varphi _{\epsilon }\psi (\sigma (x))\right\} dv_{g}d\sigma =0
\end{equation*}%
and since the functions $a$, $h$, $f$ and $\varphi _{\epsilon }$ are $G$-
invariant, the integral over $M$ does not depend on $\sigma \in G.$ Then 
\begin{equation*}
\int_{M}\left\{ \left\vert \nabla \varphi _{\epsilon }\right\vert
^{p-2}\nabla _{i}\varphi _{\epsilon }\nabla ^{i}\psi +a\left\vert \varphi
_{\epsilon }\right\vert ^{p-2}\varphi _{\epsilon }\psi -h\left\vert \varphi
_{\epsilon }\right\vert ^{q-1}\varphi _{\epsilon }\psi -f\left\vert \varphi
_{\epsilon }\right\vert ^{p^{\ast }-2-\epsilon }\varphi _{\epsilon }\psi
\right\} dv_{g}=0
\end{equation*}%
thus $\varphi _{\epsilon }$ is a weak solution of the equation(\ref{16}) in $%
W^{1,p}(M)$.

By the regularity theorem (Theorem\ref{th1}) , $\varphi _{\epsilon }\in
C^{1,\alpha }(M)\cap W^{1,p}(M)$, consequently $\varphi _{\epsilon }\mid
_{\partial M}=0$.

Now we are going to construct by mean of $\varphi _{\epsilon }$ a nodal
solution of the subcritical equation(\ref{16}). The construction is the same
as in (\cite{5}). Define%
\begin{equation*}
\psi _{\epsilon }(x)=\left\{ 
\begin{array}{c}
\left\vert \varphi _{\epsilon }\right\vert \text{ in }\Omega _{1} \\ 
-\left\vert \varphi _{\epsilon }\right\vert \text{ in }\Omega _{2}%
\end{array}%
.\right.
\end{equation*}%
Since the set $\digamma $ of fixed points of $\tau $ is negligible set, it
follows that $\psi _{\epsilon }\in H$ and arguing as above, there is $%
t_{o}>0 $ such that $J_{\epsilon }(t_{o}\psi _{\epsilon })<0$ for any $%
\epsilon \leq \epsilon _{o}$. Denote by $C^{\prime }$ the set of continuous
paths joining $0$ to $\psi _{o,\epsilon }=$ $t_{o}\psi _{\epsilon }$ and let 
\begin{equation*}
c_{\epsilon }^{\prime }=\inf_{\gamma \in C^{\prime }}\max_{t\in \left[ 0,1%
\right] }J_{\epsilon }(\gamma (t))\text{.}
\end{equation*}%
For any positive integer $m$, there exists a path $\overset{\_}{\gamma }%
_{m}\in C$ such that%
\begin{equation*}
\max_{s\in \left[ 0,1\right] }J_{\epsilon }(\overset{\_}{\gamma }%
_{m}(s))\leq c_{\varepsilon }+\frac{1}{m}\text{.}
\end{equation*}%
We let as in (\cite{5}) 
\begin{equation*}
\gamma _{m}(s)=\left\{ 
\begin{array}{c}
-\left\vert \overset{\_}{\gamma }_{m}(s)\right\vert \text{ \ \ \ in \ \ \ }%
\Omega _{1} \\ 
\ \ \left\vert \overset{\_}{\gamma }_{m}(s)\right\vert \text{ \ \ \ \ in \ \
\ }\Omega _{2}%
\end{array}%
\right.
\end{equation*}%
Clearly $\gamma _{m}\ $is a continuous path in the space $H$ and let $%
s_{m}\in \left[ 0,1\right] $ such that 
\begin{equation*}
J_{\epsilon }(\gamma _{m}(s_{m}))=\max_{s\in \left[ 0,1\right] }J_{\epsilon
}(\gamma _{m}(s))
\end{equation*}%
\begin{equation*}
\leq c_{\varepsilon }+\frac{1}{m}\text{.}
\end{equation*}
Now by the deformation lemma, there exists a continuous map 
\begin{equation*}
\eta _{m}:\left[ 0,1\right] \times H\rightarrow H
\end{equation*}%
\ such that

\begin{equation*}
\text{(i)}\ \ \ \ \ \ \ \ \ \ \eta _{m}(t,\gamma _{m}(s_{m}))=\ \gamma
_{m}(s_{m})\text{ \ for \ all \ }t\in \left[ 0,1\right]
\end{equation*}
and 
\begin{equation*}
\gamma _{m}(s_{m})\notin J_{\epsilon }^{-1}(\left[ c_{\epsilon }-\frac{1}{m}%
,c_{\epsilon }+\frac{1}{m}\right]
\end{equation*}

\begin{equation*}
\text{(ii)}\ \ \ \ \ \ 0\leq J_{\epsilon }(\gamma _{m}(s_{m}))-J_{\epsilon
}(\eta _{m}(t,\gamma _{m}(s_{m})))\leq \frac{1}{m}\text{ for all }t\in \left[
0,1\right]
\end{equation*}

\begin{equation*}
\text{(iii)\ \ \ \ }\left\Vert \ \eta _{m}(t,\gamma _{m}(s_{m}))-\ \gamma
_{m}(s_{m})\right\Vert \leq \frac{1}{m}\ \ \ \ \text{for\ \ \ \ all\ \ \ }%
t\in \left[ 0,1\right] \text{.}
\end{equation*}

\begin{equation*}
\text{(iv) \ If \ \ \ \ \ \ \ }J_{\epsilon }(\gamma _{m}(s_{m}))\leq
c_{\varepsilon }+\frac{1}{m}
\end{equation*}%
then according to the deformation lemma either 
\begin{equation*}
J_{\epsilon }(\eta _{m}(1,\gamma _{m}(s_{m})))\leq c_{\varepsilon }-\frac{1}{%
m}
\end{equation*}%
or for some \ $t_{m}\in \left[ 0,1\right] $%
\begin{equation*}
\left\Vert J_{\epsilon }^{\prime }(\eta _{m}(t_{m},\gamma
_{m}(s_{m})))\right\Vert \leq \frac{1}{m}\text{.}
\end{equation*}%
Now since we have 
\begin{equation*}
J_{\epsilon }(\gamma _{m}(s_{m}))\leq c_{\varepsilon }+\frac{1}{m}
\end{equation*}%
by (iv) we get%
\begin{equation*}
J_{\epsilon }(\eta _{m}(1,\gamma _{m}(s_{m})))\leq c_{\varepsilon }-\frac{1}{%
m}\text{.}
\end{equation*}%
Consequently since the path $s\rightarrow \eta _{m}(1,\gamma _{m}(s))$ joins 
$0$ to $\psi _{o,\epsilon }$, we obtain from the definition of $%
c_{\varepsilon }$ that

\begin{equation*}
J_{\epsilon }(\gamma _{m}(s_{m}))\geq c_{\varepsilon }\text{.}
\end{equation*}%
So the first part of (iv) cannot occur and then for some $t_{m}\in \left[ 0,1%
\right] $ 
\begin{equation*}
\left\Vert J_{\epsilon }^{\prime }(\eta _{m}(t,\gamma
_{m}(s_{m})))\right\Vert \leq \frac{1}{m}\text{.}
\end{equation*}%
Resuming, there exists $t_{m}\in \left[ 0,1\right] $ such that 
\begin{equation*}
c_{\epsilon }\leq J_{\epsilon }(\eta _{m}(t_{m},\gamma _{m}(s_{m})))\leq
c_{\epsilon }+\frac{1}{m}
\end{equation*}%
and letting 
\begin{equation*}
\varphi _{m}=\eta _{m}(t_{m},\gamma _{m}(s_{m}))
\end{equation*}%
we get a sequence of elements of $H$ \ such that 
\begin{equation*}
J_{\epsilon }(\varphi _{m})\rightarrow c_{\epsilon }\text{ \ and \ \ }%
J_{\epsilon }^{\prime }(\varphi _{m})\rightarrow 0\text{ .}
\end{equation*}%
Then as in the beginning of the proof of the Proposition\ref{p1}, \ there is
a subsequence of the sequence $(\varphi _{m})$ still denoted $(\varphi _{m})$
which converges strongly in $L^{p^{\ast }-\epsilon }(M)$ to a weak solution $%
\varphi _{\epsilon }$ of the subcritical equation. Now by (iii) $\gamma
_{m}(s_{m})\rightarrow \varphi _{\epsilon }$ strongly in $L^{p^{\ast
}-\epsilon }(M)$ then the convergence is also pointwise almost everywhere in 
$M$. Therefore $\varphi _{\epsilon }\geq 0$ on $\Omega _{1}$ and $\varphi
_{\epsilon }\leq 0$ on $\Omega _{2}$. \ Choosing a constant $B$ such that
the function $h(x,r)=a(x)\left\vert r\right\vert ^{p-1}-f(x)\left\vert
r\right\vert ^{p^{\ast }-1}-h(x)\left\vert r\right\vert ^{q}+B\left\vert
r\right\vert ^{p-1}\geq 0$ on $M\times R$ \ \ where $\left\vert r\right\vert
\leq \left\Vert \varphi _{\epsilon }\right\Vert _{L^{\infty }(M)}$ we obtain
that $\Delta _{p}\varphi _{\epsilon }+B\varphi _{\epsilon }^{p-1}\geq 0$ in $%
\Omega _{1}$ and by the strong maximum principle( Theorem\ref{th2}) we get
that $\varphi _{\epsilon }>0$ in $\Omega _{1}$, and also we have $\varphi
_{\epsilon }<0$ in $\Omega _{2}$.
\end{proof}

\subsection{The critical case}

Now we are going to show that the critical equation(\ref{1}) has a nodal
solution. First we state

\begin{proposition}
\label{p2} Let\ $G$ be a compact subgroup of the isometry group of ($M,g)$, $%
n\geq 3$, let $\tau $ be an involutive isometry of ($M,g$) such that $G$ and 
$\tau $ commute weakly and such that for some $x_{1}\in M$ \ $\tau
(O_{G}(x_{1}))\cap O_{G}(x_{1})=\phi .$ Let also $a$, $f$ and $h$ be three
smooth $G$-invariant and $\tau $-invariant \ functions. We assume that $f$
is positive on $M$ and that the operator $\varphi \rightarrow \Delta
_{p}\varphi +a\left\vert \varphi \right\vert ^{p-2}\varphi $ is coercive on $%
H$. We set $N=p^{\ast }-1$ and $q\in (p-1,N)$. \ Assume that the sequence $%
(\varphi _{\epsilon })_{\epsilon }$ of solutions of the subcritical
equations(\ref{9}) admits a subsequence which converges in $L^{k}(M)$, $k>1$%
, to $\psi \neq 0$. Then there exists $\varphi \in C^{1,\alpha }(M)$, $G$%
-invariant and $\tau $-antisymmetric in $M$ and $\varphi =0$ on $\partial M$
which is a nodal weak solution of the critical equation 
\begin{equation*}
\Delta _{p}\varphi +a\left\vert \varphi \right\vert ^{p-2}\varphi
=f\left\vert \varphi \right\vert ^{p^{\ast }-2}\varphi +h\left\vert \varphi
\right\vert ^{q-2}\varphi \text{.}
\end{equation*}
\end{proposition}

\begin{proof}
We first show that the set $(\varphi _{\epsilon })_{\epsilon }$, $\epsilon
\leq \epsilon _{o}$ of solutions to the subcritical equation(\ref{9}) is
bounded in $W^{1,p}(M)$. Let $J$ be the functional defined on the Sobolev
space $W^{1,p}(M)$ by 
\begin{equation*}
J(\varphi )=\int_{M}\left\{ \frac{1}{p}\left\vert \nabla \varphi \right\vert
^{p}+\frac{1}{p}a\left\vert \varphi \right\vert ^{p}-\frac{1}{p^{\ast }}%
f\left\vert \varphi \right\vert ^{p^{\ast }}-\frac{1}{q+1}h\left\vert
\varphi \right\vert ^{q+1}\right\} dv_{g}\text{,}
\end{equation*}%
$c=\inf_{\gamma \in C}\max_{t\in \left[ 0,1\right] }J(\gamma (t))$, where $C$
\ denotes the set of paths $\gamma $ joining $0$ and $\psi $ where $\psi $
is the function given by Lemma\ref{1}. With the same notations as in the
proof of Proposition\ref{p1}, we have%
\begin{equation*}
c_{\epsilon }=J_{\epsilon }(\varphi _{\epsilon })\leq J(\varphi _{\epsilon
})+\frac{1}{p^{\ast }}\int_{M}f\left\vert \left\vert \varphi _{\epsilon
}\right\vert ^{p^{\ast }-\epsilon }-\left\vert \varphi _{\epsilon
}\right\vert ^{p^{\ast }}\right\vert dv_{g}
\end{equation*}%
\begin{equation*}
\leq \inf_{u\in C}\max_{t\in \left[ 0,1\right] }J(u(t))+\frac{1}{p^{\ast }}%
\max_{M}f\int_{M}\left\vert \left\vert \varphi _{\epsilon }\right\vert
^{p^{\ast }-\epsilon }-\left\vert \varphi _{\epsilon }\right\vert ^{p^{\ast
}}\right\vert dv_{g}\text{.}
\end{equation*}%
So, since $\varphi _{\epsilon }\in C$, 
\begin{equation*}
c_{\epsilon }\leq c+\frac{1}{p^{\ast }}\max_{M}f\max_{t\in \left[ 0,1\right]
}\int_{M}t^{p^{\ast }-\epsilon }\left\vert \left\vert \psi \right\vert
^{p^{\ast }-\epsilon }-t^{\epsilon }\left\vert \psi \right\vert ^{p^{\ast
}}\right\vert dv_{g}\text{.}
\end{equation*}%
Then 
\begin{equation*}
\lim_{\epsilon \rightarrow 0^{+}}\sup c_{\epsilon }\leq c
\end{equation*}%
and the set $\left( \varphi _{\epsilon }\right) _{\epsilon }$ is bounded. So
there is a sequence $\left( \varphi _{n}\right) _{n}$ such that $J^{\prime
}(\varphi _{n})=0$ and $J(\varphi _{n})\rightarrow c^{\prime }.$ By
arguments as in the proof of Lemma\ref{lem2}, it follows that the sequence $%
\left( \varphi _{n}\right) _{n}$ is bounded in $W^{1,p}(M)$ and we have

$\varphi _{n}\rightarrow \varphi $ \ \ weakly in $W^{1,p}(M)$

$\varphi _{n}\rightarrow \varphi $ \ \ strongly in $L^{r}(M)$ for $r<p^{\ast
}$

$\varphi _{n}\rightarrow \varphi $ \ \ \ pointwise a.e. in $M$.

Consequently $\left\vert \varphi _{n}\right\vert ^{p^{\ast }-2}\varphi
_{n}\rightarrow \left\vert \varphi \right\vert ^{p^{\ast }-2}\varphi $
pointwise a.e. in $M$, and the sequence $\left\vert \varphi _{n}\right\vert
^{p^{\ast }-2}\varphi _{n}$ is bounded in $\left( L^{p^{\ast }}\right)
^{\prime }$ then by a well known theorem $\left\vert \varphi _{n}\right\vert
^{p^{\ast }-2}\varphi _{n}\rightarrow \left\vert \varphi \right\vert
^{p^{\ast }-2}\varphi $ weakly in $\left( L^{p^{\ast }}\right) ^{\prime }$.
The same is also true for the sequence$\left( \left\vert \varphi
_{n}\right\vert ^{q-1}\varphi _{n}\right) _{n}$ in $\left( L^{q+1}\right)
^{\prime }$ and $\varphi $ is weak solution of the critical equation(\ref{1}%
). The remaining of the proof is the same as in the second part of the proof
of Proposition\ref{p1}.
\end{proof}

To show that the sequence $(\varphi _{\epsilon })_{\epsilon }$ of solutions
of the subcritical equations(\ref{9}) admits a subsequence which converges
to $\varphi \neq 0$ in $L^{k}(M)$, $k>1$, we state.

\begin{lemma}
\label{lem3} Suppose that

(i) every subsequence of \ a sequence \ $(u_{\epsilon })_{\epsilon }$ in $%
W^{1,p}(M)$ which converges in $L^{k}(M)$ , with $k>1$, converges to $0$ \ \
\ \ \ \ \ \ \ \ \ (ii) For all $x\in M$ ,\ we can find $\delta >0$ such that%
\begin{equation}
K(n,p)^{p}(f(x))^{\frac{p}{p^{\ast }}}\underset{\epsilon \rightarrow 0^{+}}{%
\lim }\sup \left( \int_{B_{x}(\delta )\cap (B_{x}(\delta )-\partial
M)}f\left\vert u_{\epsilon }\right\vert ^{p^{\ast }-\epsilon }dv(g)\right) ^{%
\frac{p^{\ast }-p}{p^{\ast }}}<1\text{.}  \tag{17}  \label{17}
\end{equation}%
where $B_{x}(\delta )$ is the ball centred at $x$ and of radius $\delta $. \
\ \ \ \ \ \ Then \ for any $x\in M$ there is $\delta =\delta (x)>0$ \ such
that%
\begin{equation*}
\underset{\epsilon \rightarrow 0}{\lim }\sup \left( \int_{B_{x}(\delta )\cap
(B_{x}(\delta )-\partial M)}f\left\vert u_{\epsilon }\right\vert ^{p^{\ast
}-\epsilon }dv(g)\right) =0
\end{equation*}
\end{lemma}

\begin{proof}
Assume by contradiction that there is a $x_{o}\in M$ such that for any $%
\delta >0$, $\lim \sup_{\epsilon \rightarrow 0^{+}}\left(
\int_{B_{x_{o}}(\delta )\cap (B_{xo}(\delta )-\partial M)}f\left\vert
u_{\epsilon }\right\vert ^{p^{\ast }-\epsilon }dv_{g}\right) >0$. Using H%
\"{o}lder's inequality, we get 
\begin{equation*}
\int_{B_{x_{o}}(\delta )\cap (B_{xo}(\delta )-\partial M)}f\left\vert
u_{\epsilon }\right\vert ^{p^{\ast }-\epsilon }dv_{g}\leq C\left(
\int_{B_{x_{o}}(\delta )\cap (B_{xo}(\delta )-\partial M)}\left\vert
u_{\epsilon }\right\vert ^{p^{\ast }}dv_{g}\right) ^{1-\frac{\epsilon }{%
p^{\ast }}}
\end{equation*}%
where $C$ is a constant independent of $\epsilon $, and for any $s>1$ 
\begin{equation*}
\left( \int_{B_{x_{o}}(\delta )\cap (B_{xo}(\delta )-\partial M)}\left\vert
u_{\epsilon }\right\vert ^{p^{\ast }}dv_{g}\right) ^{1-\frac{\epsilon }{%
p^{\ast }}}\leq \left( \int_{B_{x_{o}}(\delta )\cap (B_{xo}(\delta
)-\partial M)}\left\vert u_{\epsilon }\right\vert ^{\frac{n(s+p-1)}{ns-p}%
}dv_{g}\right) ^{\frac{\left( p^{\ast }-\epsilon \right) (ns-p)}{p^{\ast
}(n(s+p-1))}}
\end{equation*}%
\begin{equation*}
\times \left( \int_{B_{x_{o}}(\delta )\cap (B_{xo}(\delta )-\partial
M)}\left\vert u_{\epsilon }\right\vert ^{\frac{n(s+p-1)}{n-p}}dv_{g}\right)
^{\frac{\left( p^{\ast }-\epsilon \right) \left( n(p-1)+p\right) }{p^{\ast
}(n(s+p-1))}}\text{.}
\end{equation*}%
Consequently%
\begin{equation*}
\lim_{\epsilon \rightarrow 0^{+}}\sup \int_{B_{x_{o}}(\delta )\cap
(B_{xo}(\delta )-\partial M)}\left\vert u_{\epsilon }\right\vert ^{\frac{%
n(s+p-1)}{ns-p}}dv_{g}>0
\end{equation*}%
a contradiction with the fact that any subsequence of the sequence $%
(u_{\epsilon })_{\epsilon }$ which converges in $L^{k}(M)$, for $k>1$,
converges to $0$.
\end{proof}

Now we are in position to prove Theorem\ref{th4}.

\begin{proof}
( Proof of Theorem\ref{th4}) We show that the condition$(i)$ of Lemma\ref%
{lem3} does not occur under the condition$(ii)$.

Suppose by absurd that the condition$(i)$ holds then%
\begin{equation*}
\lim_{\epsilon \rightarrow 0^{+}}\int_{M}h\left\vert u_{\epsilon
}\right\vert ^{q+1}dv_{g}=0\text{.}
\end{equation*}%
According to Lemma\ref{lem3} for every $x\in M$, there is $\delta (x)>0$
such that 
\begin{equation*}
\underset{\epsilon \rightarrow 0^{+}}{\lim }\sup \int_{B_{x}(\delta )\cap
(B_{x}(\delta )-\partial M)}f\left\vert u_{\epsilon }\right\vert ^{p^{\ast
}-\epsilon }dv_{g}=0\text{.}
\end{equation*}%
Now, since $\ M$ \ is compact, there exist $x_{1},...,x_{s}\in M$ such that $%
M=\cup _{1\leq i\leq s}B_{x_{i}}(\delta _{i}(x_{i}))$.

Consequently 
\begin{equation*}
\underset{\epsilon \rightarrow 0^{+}}{\lim }\sup c_{\epsilon }=\underset{%
\epsilon \rightarrow 0^{+}}{\lim }\sup J_{\epsilon }(u_{\epsilon })
\end{equation*}%
\begin{equation*}
=\underset{\epsilon \rightarrow 0^{+}}{\lim }\sup \int_{M}\left( \frac{q+1-p%
}{p(q+1)}h\left\vert u_{\epsilon }\right\vert ^{q+1}+\frac{p^{\ast }-p}{%
p^{\ast }}f\left\vert u_{\epsilon }\right\vert ^{p^{\ast }-\epsilon }\right)
dv_{g}=0
\end{equation*}%
which contradicts the fact that for any $\epsilon $ with $0<\epsilon \leq
\epsilon _{o}$, $c_{\epsilon }\geq \rho >0$.

So there exists $x_{o}\in M$ such that for any small $\delta >0$, 
\begin{equation}
K(n,p)^{p}(f(x_{o}))^{\frac{p}{p^{\ast }}}\lim_{\epsilon \rightarrow
0^{+}}\sup \left( \int_{B_{x}(\delta )\cap (B_{x}(\delta )-\partial
M)}f\left\vert u_{\epsilon }\right\vert ^{p^{\ast }-\epsilon }dv_{g}\right)
^{\frac{p^{\ast }-p}{p^{\ast }}}\geq 1\text{.}  \tag{18}  \label{18}
\end{equation}%
This gives%
\begin{equation*}
\lim_{\epsilon \rightarrow 0}\sup c_{\epsilon }\geq
\end{equation*}%
\begin{equation*}
\frac{p^{\ast }-p}{p^{\ast }}\lim_{\epsilon \rightarrow 0}\sup
\int_{B_{x}(\delta )\cap (B_{x}(\delta )-\partial M)}f\left\vert u_{\epsilon
}\right\vert ^{p^{\ast }-\epsilon }dv_{g}\geq \frac{p}{n}f(x_{o})^{1-\frac{n%
}{p}}K(n,p)^{-n}\text{.}
\end{equation*}

Now if $CardO_{\left\langle G,\tau \right\rangle }(x_{o})=+\infty $ we let $%
C>0$ be some given constant and we choose $\delta >0$ such that 
\begin{equation*}
C\lim_{\epsilon \rightarrow 0^{+}}\sup \int_{B_{x_{o}}(\delta )\cap
(B_{x_{o}}(\delta )-\partial M)}f\left\vert u_{\epsilon }\right\vert
^{p^{\ast }-\epsilon }dv_{g}\leq \lim_{\epsilon \rightarrow 0^{+}}\sup
c_{\epsilon }\text{.}
\end{equation*}%
Now by taking $C$ such that 
\begin{equation*}
C>K(n,p)^{n}f(x_{o})^{\frac{p}{p^{\ast }}}\lim_{\epsilon \rightarrow
0^{+}}\sup c_{\epsilon }
\end{equation*}%
we have 
\begin{equation*}
K(n,p)^{n}f(x_{o})^{\frac{p}{p^{\ast }}}\lim_{\epsilon \rightarrow
0^{+}}\sup \int_{B_{x_{o}}(\delta )\cap (B_{x_{o}}(\delta )-\partial
M)}f\left\vert u_{\epsilon }\right\vert ^{p^{\ast }-\epsilon }dv_{g}<1\text{.%
}
\end{equation*}%
Which contradicts(\ref{18}).

If $CardO_{\left\langle G,\tau \right\rangle }(x_{o})<+\infty $ we choose $%
\delta >0$ small enough such that 
\begin{equation*}
CardO_{\left\langle G,\tau \right\rangle }(x_{o})\lim_{\epsilon \rightarrow
0^{+}}\sup \int_{B_{x_{o}}(\delta )\cap (B_{x_{o}}(\delta )-\partial
M)}f\left\vert u_{\epsilon }\right\vert ^{p^{\ast }-\epsilon }dv_{g}\leq
\lim_{\epsilon \rightarrow 0^{+}}\sup c_{\epsilon }
\end{equation*}%
and taking account of (\ref{18}), we obtain%
\begin{equation*}
\lim_{\epsilon \rightarrow 0^{+}}\sup c_{\epsilon }\geq
K(n,p)^{-n}f(x_{o})^{1-\frac{n}{p}}CardO_{\left\langle G,\tau \right\rangle
}(x_{o})
\end{equation*}%
and since by construction of the sequence $(c_{\epsilon })_{\epsilon }$, 
\begin{equation*}
c\geq \lim_{\epsilon \rightarrow 0^{+}}\sup c_{\epsilon }
\end{equation*}%
it follows that%
\begin{equation*}
c\geq K(n,p)^{-n}f(x_{o})^{1-\frac{n}{p}}CardO_{\left\langle G,\tau
\right\rangle }(x_{o})\text{.}
\end{equation*}%
But this contradicts the assumption (\ref{8'}) of Theorem\ref{th4}.
\end{proof}

\section{Test functions}

Let $x_{o}$ be a point at the interior of $M$ \ such that $%
f(x_{o})=\max_{x\in M}f(x)$ and $O_{G}(x_{o})\cap \tau (O_{G}(x_{o}))=\phi $%
. Let $\psi _{\eta }$ be the radial function defined by 
\begin{equation*}
\psi _{x_{o},\eta }=\left\{ 
\begin{array}{c}
f(x_{o})^{\frac{p-n}{p^{2}}}\eta ^{\frac{n-p}{p^{2}}}(\eta +r^{\frac{p}{p-1}%
})^{1-\frac{n}{p}}C(n,p)-\mu \text{\ \ \ for \ \ }0<r\leq \delta \\ 
0\ \ \ \ \ \ \ \ \ \ \ \ \ \ \ \ \ \ \ \ \ \ \ \ \ \ \ \ \ \ \ \ \ \ \ \ \ \
\ \ \ \ \ \ \ \ \ \ \ \ \text{for\ \ \ }\ r>\delta%
\end{array}%
\right.
\end{equation*}%
where $C(n,p)=\left( n(\frac{n-p}{p-1})^{p-1}\right) ^{\frac{n-p}{p^{2}}}$, $%
\mu =f(x_{o})^{\frac{p-n}{p^{2}}}\eta ^{\frac{n-p}{p^{2}}}(\eta +\delta ^{%
\frac{p}{p-1}})^{1-\frac{n}{p}}C(n,p)$ and $\delta ,\eta $ are small
positive numbers and $r$ \ is the geodesic distance function to the point $%
x_{o}$. Suppose that $O_{G}(x_{o})=\left\{ x_{1},...,x_{n}\right\} $ and
denote 
\begin{equation*}
\overline{\psi }_{x_{o},\eta }=\sum_{i=1}^{m}\left( \psi _{x_{i},\eta -}\psi
_{\tau \left( x_{i}\right) ,\eta }\right) \text{.}
\end{equation*}%
We choose $\delta $ sufficiently small so that 
\begin{equation*}
supp(\psi _{x_{i},\eta })\cap supp(\psi _{x_{j},\eta })=\phi \text{ \ \ \ if
\ \ }i\neq j
\end{equation*}%
and 
\begin{equation*}
supp(\psi _{x_{i},\eta })\cap supp(\psi _{\tau \left( x_{j}\right) ,\eta
})=\phi \text{ \ \ for any \ \ }i\neq j.
\end{equation*}
Clearly $\overline{\psi }_{x_{o},\eta }$ is $G$-invariant and $\tau $-
antisymmetric.

At this stage, we are able to prove Theorem\ref{th3}

\begin{proof}
(Proof of Theorem\ref{th3}) Theorem\ref{th3} will be proven if the condition(%
\ref{17}) of Lemma\ref{lem3} holds and a fortiori if 
\begin{equation*}
0<c<\frac{p}{n}(\max_{M}f)^{1-\frac{n}{p}}K(n,p)^{-n}CardO_{\left\langle
G,\tau \right\rangle }(x_{o})\text{.}
\end{equation*}%
and by the definition of $c$ it suffices to show that%
\begin{equation*}
I(t\overline{\psi }_{x_{o},\eta })<\frac{p}{n}(\max_{M}f)^{1-\frac{n}{p}%
}K(n,p)^{-n}\text{\ }CardO_{\left\langle G,\tau \right\rangle }(x_{o})\text{.%
}
\end{equation*}%
Now since 
\begin{equation*}
I(\overline{\psi }_{x_{o},\eta })=cardO_{\left\langle G,\tau \right\rangle
}I(\psi _{x_{o},\eta })
\end{equation*}%
we have to show that%
\begin{equation*}
I(t\psi _{x_{o},\eta })<\frac{p}{n}(\max_{M}f)^{1-\frac{n}{p}}K(n,p)^{-n}%
\text{.}
\end{equation*}%
Put for simplicity%
\begin{equation*}
\psi _{x_{o},\eta }=\psi _{\eta }\text{.}
\end{equation*}%
The goal here is to compute the expansion in $\eta $ of $I(t\psi _{\eta })$.
Now, classical computations of $\int_{M}\left\vert \nabla \psi _{\eta
}\right\vert ^{p}dv_{g}$ give 
\begin{equation*}
\int_{M}\left\vert \nabla \psi _{\eta }\right\vert
^{p}dv_{g}=C(n,p)^{p}\left( \frac{n-p}{p-1}\right) ^{p}f(x_{o})^{1-\frac{n}{p%
}}\omega _{n-1}\frac{p-1}{p}
\end{equation*}%
\begin{equation*}
\times \left[ \int_{0}^{\infty }(1+t)^{-n}t^{n(1-\frac{1}{p})}dt-\eta
^{2\left( 1-\frac{1}{p}\right) }\frac{Scal(x_{o})}{6n}\int_{0}^{\infty
}(1+t)^{-n}t^{\left( n+2\right) -\left( 1-\frac{1}{p}\right) }dt\right]
\end{equation*}%
\begin{equation*}
+o(\eta ^{2\left( 1-\frac{1}{p}\right) })\text{.}
\end{equation*}%
We Use the following relations, for any real numbers $p$, $q$ with $p>q+1$%
\begin{equation*}
I_{p}^{q}=\int_{0}^{\infty }(1+t)^{-p}t^{q}dt=\frac{\Gamma (q+1)\Gamma
(p-q-1)}{\Gamma (p)}
\end{equation*}%
where $\Gamma $ denotes the Euler function. Such relations fulfill%
\begin{equation*}
I_{n}^{(n+2)(1-\frac{1}{p})}=\frac{\Gamma \left( (n+2)(1-\frac{1}{p}%
)+1\right) \Gamma (\frac{n-3p+2}{p})}{\Gamma \left( n(1-\frac{1}{p}%
)+1\right) \Gamma (\frac{n}{p}-1)}I_{n}^{n(1-\frac{1}{p})}
\end{equation*}%
\begin{equation*}
=a(n,p)I_{n}^{n(1-\frac{1}{p})}\text{.}
\end{equation*}
We write%
\begin{equation*}
\int_{M}\left\vert \nabla \psi _{\eta }\right\vert
^{p}dv_{g}=C(n,p)^{p}\left( \frac{n-p}{p-1}\right) ^{p}f(x_{o})^{1-\frac{n}{p%
}}\omega _{n-1}\frac{p-1}{p}
\end{equation*}%
\begin{equation}
\times \left( 1-\eta ^{^{2\left( 1-\frac{1}{p}\right) }}\frac{Scal(x_{o})}{6n%
}a(n,p)\right) I_{n}^{n(1-\frac{1}{p})}+o(\eta ^{2\left( 1-\frac{1}{p}%
\right) })\text{.}  \tag{19}  \label{19}
\end{equation}%
Now, we compute $\int_{M}a(x)\psi _{\eta }^{p}$, and get%
\begin{equation*}
\int_{M}a(x)\psi _{\eta }^{p}dv_{g}=\eta ^{p-1}C(n,p)^{p}a(x_{o})f(x_{o})^{1-%
\frac{n}{p}}\omega _{n-1}\int_{0}^{+\infty }(1+t^{\frac{p}{p-1}%
})^{p-n}t^{n-1}dt
\end{equation*}%
\begin{equation*}
+o(\eta ^{2(1-\frac{1}{p})})
\end{equation*}%
\begin{equation*}
=\frac{p-1}{p}\eta ^{p-1}C(n,p)^{p}a(x_{o})f(x_{o})^{1-\frac{n}{p}}\omega
_{n-1}I_{n-p}^{n(1-\frac{1}{p})-1}+o(\eta ^{2(1-\frac{1}{p})})\text{.}
\end{equation*}%
Taking into account the following equalities%
\begin{equation*}
I_{n-p}^{n(1-\frac{1}{p})-1}=\frac{\Gamma (n)\Gamma (\frac{n}{p}-p)}{n(1-%
\frac{1}{p})\Gamma (n-p)\Gamma (\frac{n}{p}-1)}I_{n}^{n(1-\frac{1}{p}%
)}=b(n,p)I_{n}^{n(1-\frac{1}{p})}
\end{equation*}%
we obtain%
\begin{equation*}
\int_{M}a(x)\psi _{\eta }^{p}dv_{g}=
\end{equation*}%
\begin{equation}
\frac{p-1}{p}\eta ^{p-1}C(n,p)^{p}a(x_{o})f(x_{o})^{1-\frac{n}{p}}\omega
_{n-1}b(n,p)I_{n}^{n(1-\frac{1}{p})}  \tag{20}  \label{20}
\end{equation}%
\begin{equation*}
+o(\eta ^{2(1-\frac{1}{p})})\text{.}
\end{equation*}%
Finally we compute $\int_{M}f\psi _{\eta }^{p^{\ast }}dv_{g}$ and get%
\begin{equation*}
\int_{M}f(x)\psi _{\eta }^{p^{\ast }}(x)dv_{g}=C(n,p)^{p^{\ast }}f(x_{o})^{1-%
\frac{n}{p}}\omega _{n-1}\frac{p-1}{p}\left[ I_{n}^{n(1-\frac{1}{p}%
)-1}\right.
\end{equation*}%
\begin{equation*}
\left. -\eta ^{2(1-\frac{1}{p})}\left( \frac{\Delta f(x_{o})}{2nf(x_{o})}+%
\frac{Scal(x_{o})}{6n}\right) I_{n}^{(n+2)(1-\frac{1}{p})-1}+o\left( \eta
^{2(1-\frac{1}{p})}\right) \right]
\end{equation*}%
and since%
\begin{equation*}
I_{n}^{(n+2)(1-\frac{1}{p})-1}=\frac{\Gamma ((n+2)(1-\frac{1}{p}))\Gamma (%
\frac{n+2}{p}-2)}{\Gamma (n(1-\frac{1}{p}))\Gamma (\frac{n}{p})}I_{n}^{n(1-%
\frac{1}{p})-1}
\end{equation*}%
\begin{equation*}
=c(n,p)I_{n}^{n(1-\frac{1}{p})-1}
\end{equation*}%
we can write%
\begin{equation*}
\int_{M}f(x)\psi _{\eta }^{p^{\ast }}(x)dv_{g}=C(n,p)^{p^{\ast }}f(x_{o})^{1-%
\frac{n}{p}}\omega _{n-1}\frac{p-1}{p}I_{n}^{n(1-\frac{1}{p})-1}
\end{equation*}%
\begin{equation*}
\times \left[ 1-\eta ^{2(1-\frac{1}{p})}\left( \frac{\Delta f(x_{o})}{%
2nf(x_{o})}+\frac{Scal(x_{o})}{6n}\right) c(n,p)+o\left( \eta ^{2(1-\frac{1}{%
p})}\right) \right] \text{.}
\end{equation*}%
Now letting in mind the equality 
\begin{equation*}
I_{n}^{n(1-\frac{1}{p})-1}=\frac{n-p}{n\left( p-1\right) }I_{n}^{n(1-\frac{1%
}{p})}
\end{equation*}%
we obtain%
\begin{equation*}
\int_{M}f(x)\psi _{\eta }^{p^{\ast }}(x)dv_{g}=C(n,p)^{p^{\ast }}f(x_{o})^{1-%
\frac{n}{p}}\omega _{n-1}\frac{n-p}{np}I_{n}^{n(1-\frac{1}{p})}
\end{equation*}%
\begin{equation}
\times \left[ 1-\eta ^{2(1-\frac{1}{p})}\left( \frac{\Delta f(x_{o})}{%
2nf(x_{o})}+\frac{Scal(x_{o})}{6n}\right) c(n,p)+o\left( \eta ^{2(1-\frac{1}{%
p})}\right) \right] \text{.}  \tag{21}  \label{21}
\end{equation}%
Also we have, for any real $q$ such that $\frac{n(p-1)+2p}{n-p}<q+1<p^{\ast
} $ 
\begin{equation*}
\int_{M}h(x)\psi _{\eta }^{q+1}(x)dv_{g}=\frac{p-1}{p}C(n,p)^{q+1}\eta ^{%
\frac{n-p}{p^{2}}(q+1)}f(x_{o})^{\frac{p-n}{p^{2}}(q+1)}\omega _{n-1}
\end{equation*}%
\begin{equation*}
\times \left[ h(x_{o})\eta ^{(1-\frac{n}{p})(q+1)+n(1-\frac{1}{p})}I_{(\frac{%
n}{p}-1)(q+1)}^{(n-1)(1-\frac{1}{p})-\frac{1}{p}}-\right.
\end{equation*}%
\begin{equation*}
\left. \left( \frac{\Delta h(x_{o})}{2n}+\frac{h(x_{o})Scal(x_{o})}{6n}%
\right) \eta ^{(1-\frac{n}{p})(q+1)+(n+2)(1-\frac{1}{p})}I_{(\frac{n}{p}%
-1)(q+1)}^{n(1-\frac{1}{p})+1}\right]
\end{equation*}%
\begin{equation*}
+\int_{0}^{\delta }(\eta +r^{\frac{p}{p-1}})^{(1-\frac{n}{p}%
)(q+1)}r^{n+1}dr.o(\eta ^{\frac{n-p}{p^{2}}(q+1)})+o((\eta +r^{\frac{p}{p-1}%
})^{\frac{n-p}{p}}
\end{equation*}%
\begin{equation*}
=\frac{p-1}{p}C(n,p)^{q+1}\eta ^{\frac{n-p}{p^{2}}(q+1)+(1-\frac{n}{p}%
)(q+1)+n(1-\frac{1}{p})}f(x_{o})^{\frac{p-n}{p^{2}}(q+1)}
\end{equation*}%
\begin{equation*}
\times \omega _{n-1}I_{(\frac{n}{p}-1)(q+1)}^{(n-1)(1-\frac{1}{p})-\frac{1}{p%
}}\left[ h(x_{o})-\right.
\end{equation*}%
\begin{equation*}
\left. \left( \frac{\Delta h(x_{o})}{2n}+\frac{h(x_{o})Scal(x_{o})}{6n}%
\right) \eta ^{2(1-\frac{1}{p})}e(n,p,q)\right]
\end{equation*}%
\begin{equation*}
+\int_{0}^{\delta }(\eta +r^{\frac{p}{p-1}})^{(1-\frac{n}{p}%
)(q+1)}r^{n+1}dr.o(\eta ^{\frac{n-p}{p^{2}}(q+1)})+o((\eta +r^{\frac{p}{p-1}%
})^{\frac{n-p}{p}}
\end{equation*}%
where $e(n,p,q)$ is a constant.

Since $q+1<p^{\ast }=\frac{np}{n-p}$, we have 
\begin{equation*}
\frac{n-p}{p^{2}}(q+1)+(1-\frac{n}{p})(q+1)+n(1-\frac{1}{p})
\end{equation*}%
\begin{equation*}
=\frac{np-(n-p)(q+1)}{p}(1-\frac{1}{p})>0\text{.}
\end{equation*}%
We recall that the function $u:x\in R^{n}\rightarrow C(n,p)(1+\left\vert
x\right\vert ^{\frac{p}{p-1}})^{1-\frac{n}{p}}$ realizes the equality in the
embedding $H_{1}^{p}(R^{n})\subset L_{p}(R^{n})$ that is 
\begin{equation*}
\int_{R^{n}}\left\vert \nabla u\right\vert ^{p}dx=K(n,p)^{-\frac{n}{p}%
}\int_{R^{n}}\left\vert u\right\vert ^{p^{\ast }}dx
\end{equation*}%
So from (\ref{19}), (\ref{20}) and (\ref{21}), we get%
\begin{equation*}
I(t\psi _{\eta })=\left( t^{p}-\frac{p}{p^{\ast }}t^{p^{\ast }}\right)
C(n,p)^{p}f(x_{o})^{1-\frac{n}{p}}\omega _{n-1}\left( \frac{n-p}{p-1}\right)
^{p}\frac{p-1}{p}I_{n}^{n(1-\frac{1}{p})}\frac{p}{n}
\end{equation*}%
\begin{equation*}
+F(t,n,p,\eta )
\end{equation*}%
and taking account of 
\begin{equation*}
K(n,p)^{-n}=\left( \frac{n-p}{p-1}\right) ^{p}\frac{p-1}{p}C(n,p)^{p}\omega
_{n-1}I_{n}^{n(1-\frac{1}{p})}
\end{equation*}%
we obtain that%
\begin{equation*}
I(t\psi _{\eta })=\left( t^{p}-\frac{p}{p^{\ast }}t^{p^{\ast }}\right) \frac{%
p}{n}K(n,p)^{-n}f(x_{o})^{1-\frac{n}{p}}+F(t,n,p,\eta )+H(t,n,p,\eta )
\end{equation*}%
where%
\begin{equation*}
F(t,n,p,\eta )=t^{p}a(n,p)\frac{p-1}{p}C(n,p)^{p}\left( \frac{n-p}{p-1}%
\right) ^{p}f(x_{o})^{1-\frac{n}{p}}\omega _{n-1}I_{n}^{n(1-\frac{1}{p})}
\end{equation*}%
\begin{equation*}
\times \left\{ -\frac{Scal(x_{o})}{6n}\eta ^{2(1-\frac{1}{p})}+\left( \frac{%
p-1}{n-p}\right) ^{p}a(x_{o})\frac{b(n,p)}{a(n,p)}\eta ^{p-1}\right.
\end{equation*}%
\begin{equation*}
+t^{p^{\ast }-p}\frac{n-p}{n}\left( \frac{\Delta f(x_{o})}{2nf(x_{o})}+\frac{%
Scal(x_{o})}{6n}\right) \frac{c(n,p)}{a(n,p)}\eta ^{2(1-\frac{1}{p})}
\end{equation*}%
\begin{equation*}
+o\left( \eta ^{2(1-\frac{1}{p})}\right) +o\left( \eta ^{p-1}\right)
\end{equation*}%
and%
\begin{equation*}
H(t,n,p,\eta )=-(t,t^{q+1-p}C(n,p)^{q+1-p}
\end{equation*}%
\begin{equation*}
\times \eta ^{\frac{n-p}{p^{2}}(q+1)+(1-\frac{n}{p})(q+1)+n(1-\frac{1}{p}%
)}f(x_{o})^{\frac{p-n}{p}(1-\frac{q+1}{p})}\frac{I_{(\frac{n}{p}%
-1)(q+1)}^{(n-1)(1-\frac{1}{p})-\frac{1}{p}}}{I_{n}^{n(1-\frac{1}{p})}}
\end{equation*}%
\begin{equation*}
\left. \times \left[ h(x_{o})-\left( \frac{\Delta h(x_{o})}{2n}+\frac{%
h(x_{o})Scal(x_{o})}{6n}\right) \eta ^{2(1-\frac{1}{p})}e(n,p,q)\right]
\right\} \text{.}
\end{equation*}%
Let $t_{1}\in \left[ 0,1\right] $ such that $I(t_{1}\psi _{\eta
})=\sup_{t\in \left[ 0,1\right] }I(t\psi _{\eta })$, ( $t_{1}$ is
necessarily $>0).$ Since the function $\varphi (t)=t^{p}-\frac{p}{p^{\ast }}%
t^{p^{\ast }}$ attains its maximum on the interval $\left[ 0,1\right] $ at $%
t_{o}=1$, we get 
\begin{equation*}
Sup_{t\in \left[ 0,1\right] }I(t\psi _{\eta })<K(n,p)^{-n}f(x_{o})^{1-\frac{n%
}{p}}\frac{p}{n}\left( \frac{n}{n-p}\right) ^{\frac{n}{p}}
\end{equation*}%
provided that $F(t_{1},n,p,\eta )+H(t_{1},n,p,\eta )<0.$ The assumptions $%
h(x_{o})=0$ and $\Delta h(x_{o})\leq 0$ \ give us 
\begin{equation*}
H(t_{1},n,p,\eta )\leq 0.
\end{equation*}
It remains now to show that $F(t_{1},n,p,\eta )<0$.

1) In the case\textit{\ }$1<p<2$, \ $F(t_{1},n,p,\eta )$ is equivalent to 
\begin{equation*}
G(t_{1},n,p,\eta )=
\end{equation*}%
\begin{equation*}
=t_{1}^{p}a(n,p)\frac{p-1}{p}C(n,p)^{p}f(x_{o})^{1-\frac{n}{p}}\omega
_{n-1}I_{n}^{n(1-\frac{1}{p})}a(x_{o})\frac{b(n,p)}{a(n,p)}\eta ^{p-1}
\end{equation*}%
so we must have $a(x_{o})<0$

2) In the case $p=2$, we have%
\begin{equation*}
F(t_{1},n,2,\eta )\leq t_{1}^{2}\frac{a(n,2)}{6n}C(n,2)(n-2)^{2}f(x_{o})^{1-%
\frac{n}{2}}\omega _{n-1}I_{n}^{\frac{n}{2}}
\end{equation*}%
\begin{equation*}
\times \left[ -Scal(x_{o})+\frac{6n}{\left( n-2\right) ^{2}}a(x_{o})\frac{%
b(n,2)}{a(n,2)}\right.
\end{equation*}%
\begin{equation*}
\left. +\frac{n-2}{n}\left( 3\frac{\Delta f(x_{o})}{f(x_{o})}%
+Scal(x_{o})\right) \frac{c(n,2)}{a(n,2)}\right] \eta
\end{equation*}%
\begin{equation*}
=t_{1}^{2}\frac{a(n,2)}{6n}C(n,2)(n-2)^{2}f(x_{o})^{1-\frac{n}{2}}\omega
_{n-1}I_{n}^{\frac{n}{2}}
\end{equation*}%
\begin{equation*}
\times \left[ -Scal(x_{o})+\frac{24(n-1)}{(n+2)\left( n-2\right) }%
a(x_{o})\right.
\end{equation*}%
\begin{equation*}
\left. +\frac{n-2}{n}\left( 3\frac{\Delta f(x_{o})}{f(x_{o})}%
+Scal(x_{o})\right) \frac{(n-4)n}{(n+2)\left( n-2\right) }\right] \eta
\end{equation*}%
and then the following condition must be satisfied%
\begin{equation*}
\frac{4(n-1)\text{ }}{n-2}a(x_{o})\ -Scal(x_{o})+\left( n-4\right) \frac{%
\Delta f(x_{o})}{f(x_{o})}<0\text{.}
\end{equation*}%
3) In the case $2<p<\frac{n}{2}$, to get $F(n,p,\eta )<0$, we have to assume
that 
\begin{equation*}
\left( 1-\frac{n-p}{n}\frac{c(n,p)}{a(n,p)}\right) Scal(x_{o})>\frac{3\left(
n-p\right) }{n}\frac{\Delta f(x_{o})}{f(x_{o})}\frac{c(n,p)}{a(n,p)}\text{ \ 
}
\end{equation*}%
i.e.%
\begin{equation*}
\frac{\Delta f(x_{o})}{f(x_{o})}<\frac{p}{n-3p+2}Scal(x_{o})\text{.}
\end{equation*}
\end{proof}

\section{Nonexistence results{}}

In this section we give, by mean of a Pohozaev type identity, a nonexistence
result.

\begin{proposition}
Let $n\geq 3$ and $\Omega $ be a star-shaped smooth domain of $R^{n}$ with
respect to the origin. Let $p\in (1,n)$. Suppose that $a\geq 0$, $\partial
_{r}a\geq 0$, $\partial _{r}f\leq 0$ , $h\leq 0$ , $\partial _{r}h\leq 0$
and at least one of these inequalities is strictly then the critical
equation(\ref{1}) has no nodal solution.
\end{proposition}

\begin{proof}
A Pohozaev type identity for the p-Laplacian due to Guedda and Veron \cite{7}
reads as 
\begin{equation*}
n\int_{\Omega }H(x,u)dx+\int_{\Omega }\left\langle x,\nabla
_{x}H(x,u)\right\rangle dx+(1-\frac{n}{p})\int_{\Omega }ug(x,u)dx
\end{equation*}%
\begin{equation*}
=(1-\frac{1}{p})\int_{\partial \Omega }\left\langle x,\nu \right\rangle
\left\vert \frac{\partial u}{\partial \nu }\right\vert d\sigma
\end{equation*}%
where

\begin{equation*}
g(x,u)=-a(x)\left\vert u\right\vert ^{p-2}u+f(x)\left\vert u\right\vert
^{p^{\ast }-2}u+h(x)\left\vert u\right\vert ^{q-1}u
\end{equation*}%
and 
\begin{equation*}
H(x,u)=\int_{0}^{u}g(x,s)ds\text{.}
\end{equation*}%
$\nu $ is the unit outer normal vector field to $\partial \Omega $. A direct
computation leads to the identity%
\begin{equation*}
-p\int_{\Omega }a\left\vert u\right\vert ^{p}dx+\frac{np-(q+1)(n-p)}{q+1}%
\int_{\Omega }h\left\vert u\right\vert ^{q+1}dx
\end{equation*}%
\begin{equation*}
-\int_{\Omega }\left\langle x,\nabla _{x}a\right\rangle \left\vert
u\right\vert ^{p}dx+\frac{n-p}{n}\int_{\Omega }\left\langle x,\nabla
_{x}f\right\rangle \left\vert u\right\vert ^{p^{\ast }}dx
\end{equation*}%
\begin{equation*}
+\frac{p}{q+1}\int_{\Omega }\left\langle x,\nabla _{x}h\right\rangle
\left\vert u\right\vert ^{q+1}dx=(p-1)\int_{\partial \Omega }\left\langle
x,\nu \right\rangle \left\vert \partial _{\nu }u\right\vert d\sigma
\end{equation*}%
and letting $r=\left\vert x\right\vert $, we get%
\begin{equation*}
-p\int_{\Omega }a\left\vert u\right\vert ^{p}dx+\frac{np-(q+1)(n-p)}{q+1}%
\int_{\Omega }h\left\vert u\right\vert ^{q+1}dx
\end{equation*}%
\begin{equation*}
-\int_{\Omega }r\partial _{r}h\left\vert u\right\vert ^{p}dx+\frac{n-p}{n}%
\int_{\Omega }r\partial _{r}f\left\vert u\right\vert ^{p^{\ast }}dx
\end{equation*}%
\begin{equation*}
+\frac{p}{q+1}\int_{\Omega }\left\vert u\right\vert ^{q+1}r\partial
_{r}hdx=(p-1)\int_{\partial \Omega }\left\langle x,\nu \right\rangle
\left\vert \partial _{\nu }u\right\vert d\sigma
\end{equation*}%
and the proof of the Proposition follows.
\end{proof}

\end{document}